\documentclass[a4paper,11pt]{article}
\usepackage{cancel}
\usepackage{amsmath}
\usepackage{amssymb}
\usepackage{array}
\usepackage{ragged2e}

\usepackage{soul}
\usepackage{rotating}
\usepackage{adjustbox}
\usepackage{bm}
\usepackage{longtable} 
\usepackage{arydshln}
\usepackage{adjustbox, lipsum}
\usepackage{enumerate}
\usepackage[ruled,vlined,linesnumbered]{algorithm2e}
\usepackage{rotating}
\usepackage[utf8]{inputenc}
\usepackage{amsfonts}
\usepackage{pstricks}
\usepackage{pdflscape}
\usepackage{graphicx}
\usepackage{soul}
\usepackage{afterpage}
\usepackage{graphics}

\usepackage{hyperref}

\usepackage{cleveref}

\usepackage{tikz}
\usepackage{pgfplots}
\usepackage{rotating}
\usepackage[numbers,sort]{natbib}
\usepackage{booktabs}
\usepackage{setspace}

\usepackage{authblk}
\usepackage{multirow}
\usepackage{mathtools}
\usepackage{scalerel}
\newtheorem{prop}{\bf Proposition}[section]
\newtheorem{coro}{\bf Corollary}[section]
\newtheorem{ejem}{\bf Example}[section]

\SetKwRepeat{Do}{do}{while}

\DeclarePairedDelimiter\abs{\lvert}{\rvert}%

\newcommand{\dem}{\par \noindent{\bf Proof:} }
\newcommand{\fin}{\hfill $\square$  \par \bigskip}

\newcommand{\nameM}{SL-MCFLP\xspace}
\newcommand{\nameMprime}{SL-MCFLP$^{\prime}$\xspace}
\usepackage{xcolor}
\definecolor{gr}{RGB}{0,153,0}

\usepackage{ulem}
\title{Multi-product maximal covering second-level facility location problem}
\date{}

\author{Marta Baldomero-Naranjo}
\author{Luisa I. Martínez-Merino\footnote{Corresponding author: luisa.martinez@uca.es (L.I. Martínez-Merino)}}
\author{Antonio M. Rodríguez-Chía}

\affil{\small {Departamento de Estad\'istica e Investigaci\'on Operativa, Universidad de C\'adiz}}

\allowdisplaybreaks
\setcitestyle{authoryear,open={(},close={)}}

\makeatletter
\let\oldabs\abs
\def\abs{\@ifstar{\oldabs}{\oldabs*}}
\makeatother
\begin{document}

%
\maketitle
\begin{abstract}

This paper introduces a new hierarchical facility location model with three levels: first-level facilities which manufacture different products, second-level facilities which act as warehouses and a third-level consisting of the clients who demand the products that have been manufactured in the first level and stored in the second level. In this model, called multi-product maximal covering second-level facility location problem (\nameM),  the aim is to decide the location of the second-level facilities and the products to be stored in each of them maximizing the overall clients' satisfaction with respect their coverage. 
 To deal with this model, we introduce a Mixed Integer Linear Program (MILP) which is reinforced by some families of valid inequalities. Since some of these families have an exponential number of constraints, separation algorithms are proposed. In addition, three variants of a matheuristic procedure are developed. Computational studies are included, showing the potentials and limits of the formulation and the effectiveness of the heuristic.  

\end{abstract}

\textbf{Keywords:} hierarchical location, covering, integer programming, matheuristic.  

\section{Introduction}

Health care systems, telecommunication networks, disaster management systems, or production-distribution systems are real situations where hierarchical facility location problems (HFLP) have been successfully applied, see for example \cite{afshari2014challenges, CHAUHAN20191, paul2017multiobjective,li2018cooperative,Contreras2019} and references therein. 
HFLP aims to locate facilities that interact in a multi-level hierarchical framework, see \cite{ORTIZASTORQUIZA2018791} for a comprehensive review on multi-level facility location problems.
Particularly, these problems can be applied to an increasingly relevant activity: e-commerce, in which efficient hierarchical management of facilities (product sources, warehouses, lockers, etc.) is needed.

 Depending on the objective of the problem, these have been divided into three main types: median, covering, and fixed charge objectives, see \cite{reviewHierarchical14} for a detailed description. Among the models related to covering objectives, the main studied ones correspond with i) the classical set covering, that is, minimizing the number of open facilities ensuring a complete covering of the demand, presented in  \cite{TorSwaRevBer71}, and ii) the classical maximal covering, which consists of maximizing the covered demand considering that a limited number of facilities can be operating, 
 introduced in \cite{ChuRev74}. Since the introduction of these covering objectives five decades ago, many extensions and variants of these problems have been studied. This is mainly due to their great applicability in numerous contexts, as reflected in the wide variety of recent studies in the literature, e.g., 
 \cite{BalKalMarRod22,BalKalRod20, MarMarRodSal18,BlaGazSad23,PelXu23,VatJay21,GarMar19,FroMaiHam20}.  
 The problem studied in this paper is focused on the maximal covering on a multi-level location framework.


In this context, we introduce a new model, the multi-product maximal covering second-level facility location problem (\nameM), in which we assume a multi-product hierarchical system in a single flow-pattern  with costumer preferences. This problem considers that the first-level facilities (factories, farms, product sources, etc.) are already established and it focuses on determining the optimal location for the second-level facilities (warehouses, shops, delivery centers, etc.) which cover the demand of the clients. 

The location of second-level facilities such as warehouses or delivery centers has been previously addressed in the literature in different specific applications, see for instance  \cite{Naji2012}, \cite{HUANG2015} or \cite{MILLSTEIN2022}. Namely, \cite{Naji2012} address the location of satellite distribution centers to deliver humanitarian aid, which have to be supplied from a central depot. They present a generalization of the covering tour problem. \cite{HUANG2015} tackle the site selection and space determination for warehouses in a two-stage network (products are shipped from part suppliers to warehouses and then delivered to assembly plants). The objective is to minimize transportation and operation costs. \cite{MILLSTEIN2022} discuss a model to find the optimal locations and capacities of omnichannel warehouses, and the product flows in the system (from suppliers to warehouses, warehouses to the stores, and warehouses and stores to online markets).

In \nameM, the facilities interact in a single-flow pattern. This means that clients' demands can only be satisfied by second-level facilities within a coverage radius and, at the same time, the products offered by second-level facilities are obtained from first-level facilities within another coverage radius. That is, there is double coverage, the second-level facilities (e.g. warehouses) should be within the coverage radius ($R$) from the first-level facilities (e.g. factories) and the clients should be within the coverage radius ($r$) from the second-level facilities. Again, this assumption can be applied in e-commerce, where customers have to be served within a specified period of time and the product should be delivered from the factory to the warehouse, and later, from there to the customer. This requirement can be modeled using these double coverage radii.  

Other models where multiple coverages are  considered can be found in \cite{Landete2019} or \cite{Y2023}. Particularly, in \cite{Landete2019}, four models where cooperative users act as intermediate facilities are proposed related to humanitarian logistics. In one of these models, cooperative users present a double coverage framework: one to be served and one to serve other users. 
Besides, \cite{Y2023} deal with a problem related to the relocation of emergency service vehicles where demand points have different requirements of vehicles to be considered as covered.

In facility location problems, customer preferences are also an important feature to take into account. Usually, the quality of a product or its price depends on its  source. Consequently, the customers will vary their satisfaction with a product regarding its manufacturer (origin). For instance, \cite{Casas2020} introduce a bi-level model based on maximal covering location problem where clients are allocated to the facilities depending on their preferences. In \cite{DOLAMA2021}, an extension of the rank pricing problem is proposed. Note that, in this problem, customers have a ranked list of preferred products and a budget. Another model dealing with customer preferences is the one introduced in \cite{Mrkela2022} where customers can choose among the facilities operating in a certain radius.


 To handle clients' preferences, in the \nameM, it is assumed that the same second-level facility can offer: i) several products and ii) the same product from different origins. Moreover, each client has certain preferences for the same product depending on the first-level facility which produces that product. These preferences are reflected in the model by different weights in the objective function depending on the client and the first-level facility. Therefore, if the customer can be served by  the same product from several origins, he/she will choose the most preferred one.

 As far as we are aware, the \nameM is the first model for the location of second-level facilities included in a multi-product and single-flow covering context which also includes customers' preferences to evaluate the coverage of the clients. For an extensive review on covering models, the reader is referred to \cite{ChurchMurray2018}. 

As mentioned before, each second-level facility can offer several products. \nameM selects which products are offered in each candidate facility taking into account that the opening cost of the facility will depend on the number of different products that it can offer. In other words, these costs could be interpreted as infrastructure costs for offering each type of product,  e.g. a storage shelf, a refrigerated area, a frozen product area, among others.  The overall cost of the facilities establishment is limited in the model by a budget constraint. In addition, the number of different products offered by a second-level facility is bounded.

Therefore, the main assumptions of this new model are: 
\begin{itemize}
\item We consider a multi-product and single-flow hierarchical context, where the location of first-level facilities (sources) is given. 
There are three main decisions: i) where to locate the second-level facilities (warehouses), ii) their size (number of different products that they offer), and iii) which products offer each second-level facility.
\item There are several sources for the same product. The satisfaction of each product is different for each client and depends on its source. 
\item A client demands several products and the same product will be served by at most one source. In order to satisfy the customers' demand, there must be a double coverage: the client must be covered by a second-level facility, and this, in turn, by a first-level facility. The coverage radii depend on the type of product.
\item  A second-level facility can offer several products, although the number of different products is limited. Installing each facility has a cost depending on the number of different products offered. 
\end{itemize}

\begin{figure}[htb]\label{fig:modelo}
\centering
\includegraphics[width=14cm]{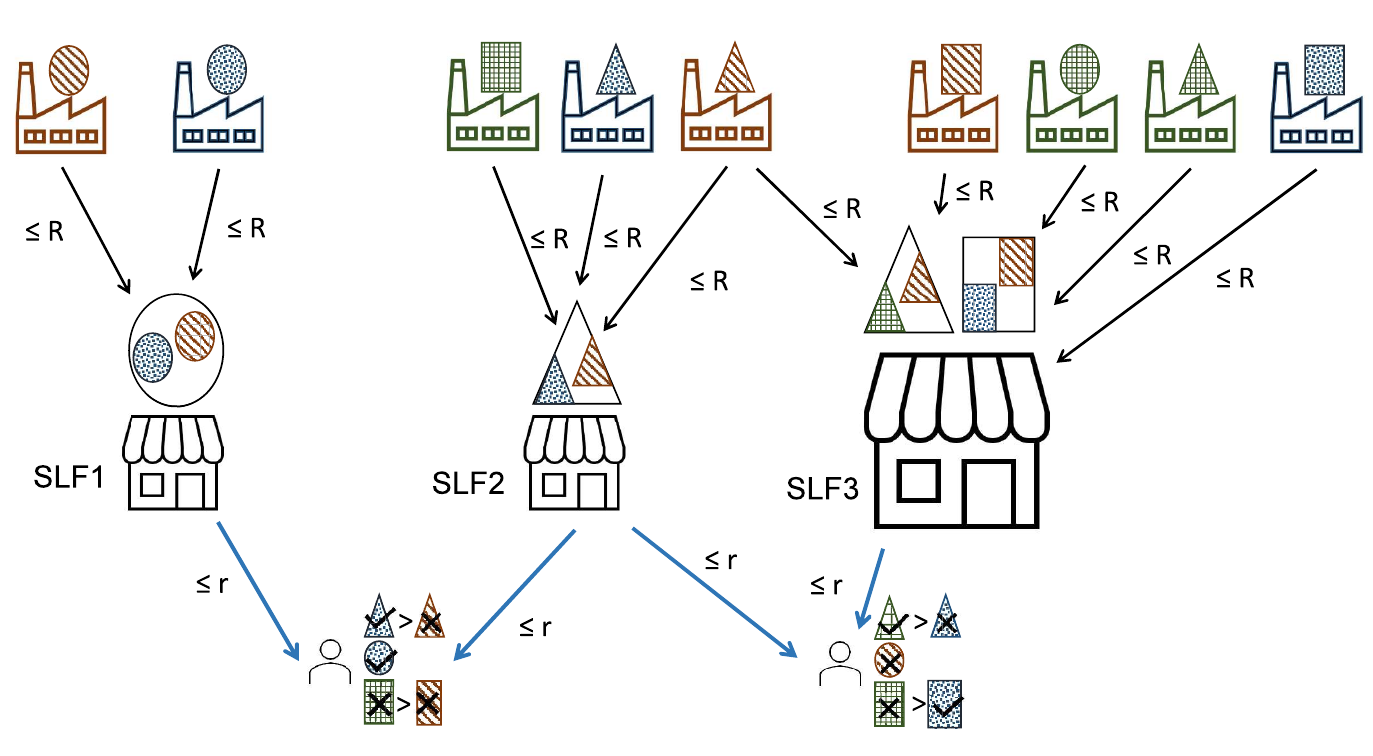}
\caption{Graphical representation of the problem.}
\end{figure}

Figure~\ref{fig:modelo}  shows a graphical representation of the problem. The factories represent the first-level facilities, where the geometrical shapes represent different products (rectangles, triangles, and ovals). The color and type of filling of the geometric shapes represent the origin (sources) of the product. The warehouses represent the second-level facilities, 
 denoted as SLF1, SLF2,  and SLF3. Their sizes in the figure are proportional to the number of different products that can offer. In fact, we assume that warehouses SLF1 and SLF2 can offer one product, and warehouse SLF3 can offer two different products. The figure shows the types of products that can be supplied in these facilities, together with the sources of the products (the latter depends on the coverage radius of the first-level facilities per product). Finally, the clients have been depicted, with their demand for the products, ordered from left to right from the highest to the lowest preference depending on the source. The tick represents the product and source served to the client and the cross represents the ones that are not served. Note that for each product, the client can only be served by one source and that clients cannot be directly served by a first-level service. Furthermore, to simplify the figure, coverage radii have been assumed to be the same for all products and only the arcs within the corresponding coverage radius have been represented. Observe that the green square shaped product is not served to the clients. Even if the coverage constraints are met, (i.e., there is a first-level facility within a distance $R$ of a second-level facility and this second-level facility is within a distance $r$ of the client), the limitation on the number of different products that each warehouse can offer  prevents the second-level facility SLF2 from offering square shaped products.  

The main contributions of this paper are the following:
\begin{itemize}
\item A new hierarchical facility location model is introduced (\nameM). Furthermore, a mixed integer linear formulation which can be solved by off-the-shelf solvers is proposed for the \nameM.
\item Some families of valid inequalities are developed to reduce the solution time of the formulation. In some cases, the number of valid inequalities is exponential and, consequently, separation  algorithms are developed along the paper.
\item Three different variants of a matheuristic procedure are introduced. This matheuristic is based on the Adaptive Kernel Search  (see \cite{AKS}) which obtains good quality solutions by solving restricted problems derived from the original one.  
\end{itemize}

The remainder of the paper is organized as follows. In Section~\ref{sec:model}, the notation is introduced and a formulation is proposed to model the problem under study. In Section~\ref{sec:VI}, some families of valid inequalities are presented for the formulation. In cases in which the number of valid inequalities is exponential, different separation procedures have been developed. In Section~\ref{sec:Heuristic}, three variants of a matheuristic algorithm are presented. In Section~\ref{sec:Computational}, an extensive computational analysis is carried out. This illustrates the usefulness of the valid inequalities and the corresponding separation methods, as well as, shows the good performance of the proposed heuristic algorithms. Then, our conclusions are presented in Section~\ref{sec:conclusion}. Lastly, a brief appendix with an alternative formulation is provided.




\section{Model}\label{sec:model}

In the \nameM, we assume that there exists a set of clients, $I$, demanding some products  manufactured by the set of first-level facilities, denoted by $S$. Each first-level facility, $s\in S,$ produces or offers one product of the set of products $M$, and this product is denoted as $m(s)$.  Besides, each client $i\in I$ has a different preference/satisfaction for a product depending on its origin $s\in S$, $w_{is}$. 

\begin{table}[htb]
    \centering
      
\begin{tabular}{lp{0.85\textwidth}}
\hline
$B$& Budget.  \\
$I$ & Set of clients.\\
$J$ & Set of potential second-level facilities.\\
$M$ & Set of products.\\
$S$ & Set of established first-level facilities.\\
$b_j$&Maximum number of different products that $j\in J$ can offer.\\
$T(j)$ & Set of numbers of different products that a second-level facility $j\in J$ can offer, i.e., $T(j)=\{0,1,\dots,b_j\}.$\\
$J_T$ & The set of pairs $(j,t)$, such that, $j \in J$ and $t \in T(j)$. \\
$m(s)$& Product that offers the first-level facility $s\in S$.\\
$d_{kl}$&Distance between $k$ and $l$ (where $k$ and $l$ can refer to clients, first and/or second-level facilities) \\
$R_m$&  Covering radius of product $m\in M$ from a first-level to a second-level facility. Thus, $j\in J$ can offer product $m(s)$ if $d_{sj}\leq R_{m(s)}$. \\
$r_m$&  Covering radius of product $m$ from a second-level facility to a client. Thus, the demand of product $m\in M$ at client $i\in I$ can be covered by $j\in J$ if $d_{ji}\leq r_{m}$. \\
$J(i,s)$& Set of second-level facilities $j\in J$ such that $d_{ji}\leq r_{m(s)}$  and  $d_{sj}\leq R_{m(s)}$.\\
$M_I(i)$& Set of products $m\in M$ that client $i\in I$ demands. \\
$M_J(j)$& Set of products $m\in M$, that can be offered from $j\in J$, i.e., those such that exists $s\in S$ with  $m(s)=m$ and $d_{sj}\leq R_{m(s)}$ \\
$S_I(i)$ & Set of $s\in S$ whose products are demanded by client $i\in I$ and such that exists $j\in J$ with  $d_{ji}\leq r_{m(s)}$ and  $d_{sj}\leq R_{m(s)}$. \\
$S_J(j)$& Set of $s\in S$ such that $d_{sj}\leq R_{m(s)}$, i.e., $m(s) \in M_J(j)$. \\
$S_M(m)$ & Set of $s\in S$ that serves product $m\in M$, i.e., $m(s)=m$. \\ 

$w_{is}$& Satisfaction of serving product $m(s)$ from $s\in S$ to client $i\in I$.\\ 
$c_{jt}$& Cost of offering $t\in T(j)$ different products in $j\in J$.\\


\hline
\end{tabular}

 \caption{Notation used in the paper for the \nameM model}
    \label{tab:notation}
\end{table}

 In this problem, we have a given set of candidate locations, $J$, for the second-level facilities. The objective of \nameM is to determine which of these facilities are opened and which products are offered by them maximizing the sum of the clients' satisfaction. 
 Let $j\in J$ be an operating second-level facility, $T(j)$  is the set of numbers of different products that $j$ can offer and $J_T:=\{(j,t) \: : \: j \in J \mbox{ and } t \in T(j)\}$. Then, $c_{jt}$ for $(j,t)\in J_T$ is the cost of setting a service facility at $j$ offering $t$ products. We assume that these costs are non-decreasing, i.e., given a second-level facility $j\in J$, the cost of offering one product ($c_{j1}$) is always lower than or equal to the cost of offering two products ($c_{j2}$), and so on. 
 There exists a budget $B$ for the maximum cost of establishing the second-level facilities.
 
 In addition, a candidate facility $j\in J$ can provide a product from a first-level facility $s\in S$ to a client $i\in I$ if and only if  $d_{sj}\leq R_{m(s)}$, the facility $j$ offers product $m(s)$, and  $d_{ji}\leq r_{m(s)}$. 
Table \ref{tab:notation} reports a detailed description of the previously presented notation and some extra notation that would be used along the paper.

To model the \nameM, we include the following families of binary variables:
\begin{eqnarray*}
y_{is}&=&\begin{cases}
1,&\mbox{if first-level facility $s$ serves} \\
&\mbox{product $m(s)$ to client $i$, }\\
0,&\mbox{otherwise,}
\end{cases} \mbox{ for } i\in I, s \in S_I(i),\\
x_{jm}&=&\begin{cases}
1,&\mbox{if second-level facility $j$ }\\
&\mbox{offers product $m$,}\\
0,&\mbox{otherwise,}
\end{cases} \hspace{0.5cm}\mbox{ for } j\in J, m\in M_J(j),\\
u_{jt}&=&\begin{cases}
1,& \mbox{if second-level facility $j$}\\
&\mbox{offers $t$ products, }\\
0,&\mbox{otherwise,} \\
\end{cases} \hspace{0.6cm}\mbox{ for } j\in J, t\in T(j).
\end{eqnarray*}

The resulting formulation is the following,

\begin{align}
\mbox{(\nameM)} \quad &\max&&\displaystyle \sum_{i\in I}\sum_{s\in S_I(i)} w_{is}y_{is}\nonumber\\
&\mbox{s.t.}&&\sum_{j\in J(i,s)}x_{jm(s)}\geq y_{is},&& i\in I, s\in S_I(i), \label{con:cubri}\\
&&&\sum_{s\in S_I(i)\cap S_M(m)}y_{is}\leq 1,&& i\in I,m\in M_I(i),\label{con:onesource}\\
&&&\sum_{m\in M_J(j)}x_{jm}=\sum_{t\in T(j)}t u_{jt},&&j\in J,\label{con:nprod}\\
&&&\sum_{t\in T(j)}u_{jt}= 1,&&j\in J,\label{con:nprod2}\\
&&&\sum_{j\in J}\sum_{t\in T(j)}c_{jt} u_{jt}\leq B,\label{con:budget}\\
&&&y_{is}\in \{0,1\},&&i\in I,s\in S_I(i),\label{con:y_b}\\
&&&u_{jt}\in \{0,1\},&&j\in J, t\in T(j), \label{con:sumut}\\
&&&x_{jm}\in \{0,1\},&&j\in J,m\in M_J(j). \label{con:x}
\end{align}
In this formulation, the objective function represents the overall satisfactions from the covered products of each client. Constraints \eqref{con:cubri} ensure that client $i$ is only covered by a certain first-level facility $s\in S$ if there is at least one second-level facility $j$ which offers $m(s)$, $d_{sj}\leq R_{m(s)}$, and $d_{ji}\leq r_{m(s)}$. 
Constraints \eqref{con:onesource} 
ensure that each customer is served by only one source per product, although he/she can demand several different products. Clearly, the customer will 
choose the source, among all feasible ones, for which he/she has the highest preference. Besides, constraints \eqref{con:nprod} and \eqref{con:nprod2} determine the number of different products that offers each second-level facility $j\in J$. Budget constraint \eqref{con:budget} determines the maximum cost of opening the second-level facilities. Finally, constraints \eqref{con:y_b} to \eqref{con:x} determine the domain of definition of the variables appearing in the model. 

Next, we introduce a result proving that the integrality condition of $y$-variables can be relaxed. 
\begin{prop}
The integrality condition of $y$-variables in constraints \eqref{con:y_b} can be relaxed. Thus, constraints \eqref{con:y_b} could be replaced by 
		\begin{equation}
		0\leq y_{is}\leq 1,    \quad   i\in I,s\in S_I(i).  \label{con:y} 
		\end{equation} 
\end{prop}
\dem
Observe that the objective function pushes $y$-variables to be as large as possible and $y$-variables are bounded by the partial sum of $x$-variables in constraints~\eqref{con:cubri}. Therefore, the integrality condition of $y$-variables is inherited from the one of $x$-variables.
If there are ties in the vector $w_{is}$, for $s\in S_I(i)$, a fractional optimal solution could be obtained for $y$-variables. In this case, an integer solution with the same objective value can be easily constructed from a non-integer optimal solution. \mbox{ }\fin



In the next section, we describe some valid inequalities for formulation \nameM which reinforce it and improve its computational performance.

\section{Valid inequalities}\label{sec:VI}

This section is focused on the development of valid inequalities. First of all, we present constraints that relate $x$- and $y$- variables. Then, valid inequalities using the budget constraint are proposed. However, some of the previously mentioned valid inequalities have exponential cardinality. Therefore, we present a procedure and an auxiliary mixed integer linear program (MILP) to identify the valid inequalities of those families that cut the current fractional solution in the branch and bound tree, i.e., we develop separation methods.

\subsection{Valid inequalities linking $x$- and $y$-variables}

In the following result, we present two families of valid inequalities. The first one establishes the relationship between i) the covering variables, $y$-variables, of a pair of clients that demands the same product and ii) the $x$-variable of the second-level facilities offering this product. The second family of constraints links, for a given client and product, the covering variables and the $x$-variable of this product for each second-level facility in the coverage area.

\begin{prop}
The following families of constraints are valid inequalities for the \nameM: 

\begin{enumerate}[i)]
    \item For each $i,i^\prime \in I$ and $s\in S_I(i)$ such that $J(i,s)\cap J(i',s)\neq \emptyset$:
\begin{equation} \label{ec:DDVVcubri}
    y_{is} \leq \sum_{s^{\prime}\in S_I(i^\prime) \cap S_M(m(s)):   w_{i^{\prime}s^{\prime}}\geq w_{i^\prime s}}y_{i^\prime s^{\prime}} + \sum_{j\in J(i,s) \setminus J(i^\prime,s) } {x}_{jm(s)}.
\end{equation}


\item   For each $i\in I, s\in S_I(i), j\in J(i,s)$:
\begin{equation}\label{ec:DDVVpref}
  x_{jm(s)} + \sum_{s^{\prime}\in S_I(i) \cap S_M(m(s)): w_{is}> w_{is^\prime}}y_{is^{\prime}} \leq 1.
\end{equation}

\end{enumerate}
\end{prop}
\dem 

\textit{i)} 
Let $i,i'\in I$ and $s\in S_I(i)$ such that $J(i,s)\cap J(i',s)\neq \emptyset$. Assume that $y_{is}=1$, then there exists $j_0\in J(i,s)$ such that $x_{j_0 m(s)}=1$, $d_{j_0 i}\leq r_{m(s)}$, and $d_{sj_0}\leq R_{m(s)}$.

If $j_0\in J(i,s)\cap J(i',s)$, then $ \displaystyle\sum_{s^{\prime}\in S_I(i^\prime) \cap S_M(m(s)):   w_{i^{\prime}s^{\prime}}\geq w_{i^\prime s}}y_{i^\prime s^{\prime}}=1$, since $d_{j_0i'}\leq r_{m(s)}$, $d_{s j_0}\leq R_{m(s)}$ and, consequently, $i^{\prime}$ will receive product $m(s)$ at least with satisfaction $w_{i's}$. On the other side, if 
$j_0\in J(i,s)\setminus J(i',s)$, then $\sum_{j\in J(i,s) \setminus J(i^\prime,s) } {x}_{jm(s)} \ge x_{j_0 m(s)}=1$.

In both cases, if $y_{is}=1$, then the right-hand side is at least one. These constraints ensure that if a client $i$ is covered by $s$,  there are two options: the client $i^\prime$ is also covered (by $s$ or by another first-level facility whose satisfaction is greater), or the client $i$ is covered by a second-level facility outside the coverage radius of $i^\prime$.

\textit{ii)} 
These constraints establish that if the second-level facility $j$ offers a product $m$ of source $s$ and client $i$ can be covered by $j$, then the demand of product $m$ of the client $i$ can not be satisfied by a first-level facility whose satisfaction is lower than the one of $s$ and vice versa. 

\fin

In the next result, we introduce a new family of valid inequalities for \nameM. Particularly, for a given client and a product, these constraints consider subsets of origins that could provide this product for this client.


\begin{prop}
The following family of constraints 
    are valid inequalities for formulation \nameM:
    \begin{equation}\label{ec:DDVVnew4}
        \sum_{s\in S_1} y_{is}\leq \sum_{j\in \underset{s\in S_1}{\cup}J(i,s)}{x}_{jm}, \quad i\in I, m\in M_I(i), S_1\subseteq S_I(i)\cap S_M(m).
    \end{equation}
\end{prop}

\dem
 Given $i\in I$, $m\in M_I(i)$ and $S_1\subseteq S_{I}(i)\cap S_M(m)$, assume that $y_{is_0}=1$ for $s_0\in S_1$. Then, due to constraints \eqref{con:cubri}, there exists $j\in J(i,s_0)$ such that $x_{jm(s_0)}=1$, i.e., $\displaystyle\sum_{j\in \underset{s\in S_1}{\cup}J(i,s)}{x}_{jm}\geq 1$. Therefore, using \eqref{con:onesource}, the result follows.
 
\fin

Note that the cardinality of the family of valid inequalities~\eqref{ec:DDVVnew4}  is exponential. Therefore, it is worth studying which subsets $S_1$ are appropriate to be taken into account.
In the following, we present two particular families of valid inequalities~\eqref{ec:DDVVnew4} whose cardinality is polynomial and that report good computational results, see Section~\ref{sec:Computational} for more details.
\begin{enumerate}[i)]
   \item  
    \begin{equation}
    \label{ec:DDVVnew1}
        \sum_{s\in S_I(i)\cap S_M(m)} y_{is}\leq \sum_{j\in \underset{\stackrel{s\in S:}{m(s)=m}}\cup J(i,s)}{x}_{jm}, \quad i\in I, m\in M_I(i).
    \end{equation}
        \item \begin{eqnarray}\label{ec:DDVVnew2}\sum_{\stackrel{s^\prime \in S_I(i)\cap S_M(m):}{ w_{is^\prime}\geq w_{is}}} \hspace{-0.8cm} y_{is^{\prime}}\leq \sum_{j\in \underset{\stackrel{s^\prime\in S_I(i)\cap S_M(m):}{ w_{is^\prime}\geq w_{is}}}{\cup}  J(i,s^\prime)}\hspace{-1cm} {x}_{jm}, \;\;\; i\in I, m\in M_I(i), s\in S_I(i)\cap S_M(m).\end{eqnarray}   
\end{enumerate}


In addition, alternative valid inequalities of the family \eqref{ec:DDVVnew4} could be derived  from fractional solutions obtained in each node  of B\&B tree. Algorithm \ref{strategy1} describes a separation procedure. For each fractional solution $(\bar{x},\bar{y},\bar{u})$ in each node of the B\&B tree of the \nameM, we build a set $S_1$ composed by the indices of the four most promising first-level facilities (lines 2-3 of Algorithm \ref{strategy1}). Note that if, for a given client and product, that customer can reach less than four first-level facilities that offer that product, the algorithm would be adjusted accordingly. If constraints~\eqref{ec:DDVVnew4} using this set $S_1$ does not hold for the current solution then, this valid inequality is included in the model (lines 4 and 5 of Algorithm \ref{strategy1}).

\begin{algorithm}[htbp]\label{strategy1}
\setstretch{1.15}
	\KwData{A fractional solution $(\bar{x},\bar{y},\bar{u})$ in the B\&B tree of \nameM.}
	\KwResult{A set of valid inequalities of family \eqref{ec:DDVVnew4}.}
\For{ $i\in I$ and $m\in M_I(i)$} { Obtain the following values:
\begin{eqnarray*}\nonumber
    s^0_{im}&:=&\underset{s\in S(m)}{\mathrm{argmax}}\{\bar{y}_{is}\},\\
    s^1_{im}&:=&\underset{s\in S(m)\setminus \{s^0_{im}\}}{\mathrm{argmax}}\left\{\bar{y}_{is}-\sum_{j\in J(i,s)\setminus J(i,s^0_{im})} \bar{x}_{jm}\right\},\\
    s^2_{im}&:=&\underset{s\in S(m)\setminus \{
    s^0_{im},s^1_{im}\}}{\mathrm{argmax}}\left\{\bar{y}_{is}-\sum_{j\in J(i,s)\setminus (J(i,s^0_{im})\cup J(i,s^1_{im}))} \bar{x}_{jm}\right\}, \\
     s^3_{im}&:=&\underset{s\in S(m)\setminus \{s^0_{im},s^1_{im},s^2_{im}\}}{\mathrm{argmax}}\left\{\bar{y}_{is}-\sum_{j\in J(i,s)\setminus (J(i,s^0_{im})\cup J(i,s^1_{im}) \cup J(i,s^2_{im}))} \bar{x}_{jm}\right\}.
    \end{eqnarray*}
The $s$-values will be discarded if the maximum values resulting from the above expressions are lower than or equal to zero. \\
Define $S^1_{im}:=\{s^0_{im},s^1_{im},s^2_{im},s^3_{im}\}$
(most promising first-level facilities).
\\
\If{$\displaystyle\sum_{s\in S^1_{im}} \bar{y}_{is}>\sum_{j\in \underset{s\in S^1_{im}}{\cup}J(i,s)}\bar{x}_{jm}$}{
Add the next valid inequality of family \eqref{ec:DDVVnew4} to the model:
$$\sum_{s\in S^1_{im}} y_{is}\leq \sum_{j\in \underset{s\in S^1_{im}}{\cup}J(i,s)}{x}_{jm}.$$
}
}
\caption{Separation procedure of family \eqref{ec:DDVVnew4}.}
\end{algorithm}

\subsection{Covering valid inequalities}
 Next, we provide a family of valid inequalities based on the covering properties of the model.



\begin{prop}\label{prop:DDVVcovMax}

The following family of constraints is valid for the formulation\\ \nameM: 
\begin{equation}
\label{con:DDVVBRef}
 \sum_{(j,t) \in J_T^1 \cup J_T^2} u_{jt} \leq  \tilde{k},
 \quad  J_T' \subseteq J_T,
 \end{equation}
such that,
\begin{eqnarray*}
J_T^1 &:=& \{(j,t) \in J_T \: : \: \exists t'\le t : (j,t') \in J_T'\}, \\
J_T^2 &:=& \{(j,t) \in J_T\setminus J_T^1 \: : \: c_{jt} > \tilde{c}\}. 
\end{eqnarray*}
where $\tilde{k}$ is the largest index $k$ that satisfies the following condition:
\begin{equation}\label{con:DDVVB}
    \sum_{h=1}^{k} \delta_{(h)}\leq B,
\end{equation}
and $\delta(\cdot)$ is a permutation of costs $c_{jt}$ for 
$(j,t) \in J_T'$
sorted in non-decreasing order and $\tilde{c}=\max\left\{B-\sum_{h=1}^{\tilde{k}}\delta_{(h)},\delta_{(\tilde{k})}\right\}$.
\end{prop}

  For the sake of clarity, we introduce an example illustrating Proposition \ref{prop:DDVVcovMax} before its proof.

\begin{ejem}\label{ex:1}
Consider five candidate second-level facilities and five products, where 
the cost $c_{jt}$ for $j=1,\dots,5$ and $t=1,\dots,5$ are given in Table~\ref{tb:ej}. 
Furthermore, the budget, B, is equal to 20. For simplicity, we consider that $b_{j}=5$, for $j=1,\dots,5.$

We set $J_T'=\{(1,4),(2,4),(3,4),(4,4)\}$.
Next, we sort the cost of $c_{jt}$ for $(j,t)\in J_T'$
in non-decreasing order, i.e, 6,7,9,9. Then, we find the largest index $k$ that satisfies condition \eqref{con:DDVVB} $(\tilde{k}=2)$ and $\tilde{c}=\max\{20-13,7\}=7$. Thus, the sets defined in the proposition would be \begin{eqnarray*}
    &&J_T^1=\{ (1,4), (1,5), (2,4), (2,5),  (3,4), (3,5), (4,4), (4,5)\}, \\
    &&J_T^2=\{(1,2), (1,3), (3,3), (5,5)\}.
\end{eqnarray*}
\begin{table}[ht]
\centering
\begin{tabular}{c|c|c|c|c|c}
     Warehouse & $t=1$ & $t=2$ & $t=3$ &$t=4$&$t=5$  \\
     \hline
     $j=1$&4&\textbf{8}&  \textbf{8.5}& \textbf{9}&\textbf{10}\\
     $j=2$&2&3&6& \textbf{7}&\textbf{9}\\
     $j=3$&1&3&\textbf{8}&\textbf{9}&\textbf{10}\\
     $j=4$&3&4&5&\textbf{6}&\textbf{7}\\
     $j=5$&1&3&4&6.5&\textbf{7.5}\\
     \hline
\end{tabular}
\caption{Cost of offering $t$ products in facility $j$ of Example \ref{ex:1}.}
\label{tb:ej}
\end{table}
 Therefore, for this example and this choice of subset
$J_T' \subseteq J_T$ the valid inequality~\eqref{con:example} is obtained from constraints \eqref{con:DDVVBRef}. Note that the $u$-variables that are included in  \eqref{con:example} are the ones related to the highlighted costs in Table~\ref{tb:ej}, i.e.,
\begin{equation}\label{con:example}
u_{12}+u_{13}+u_{14}+u_{15}+u_{24}+u_{25}+u_{33}+u_{34}+u_{35}+u_{44}+u_{45}+u_{55}\leq 2.
\end{equation}

\end{ejem}

\dem
For any given $J^\prime_T\subseteq J_T$, we select the cheapest elements of this set such that their sum is less than the budget, i.e., the maximum number of elements that satisfy condition \eqref{con:DDVVB}. In this case, by construction, $\tilde{k},$ is the maximum number of $u_{jt}$ variables for $(j,t)\in J^\prime_T$ that can take value one simultaneously in any feasible solution of the problem. Clearly, if $\tilde{k}$+1 variables of the set $J^\prime_T$ take value one, constraint \eqref{con:budget} is not fulfilled. Therefore, the following constraint is valid:
\begin{equation}\label{ec:profProp0}\sum_{(j,t) \in J_T'} u_{jt} \leq  \tilde{k}.
\end{equation}
However, the latter constraint can be reinforced. Indeed, using that the setting costs of any facility are non-decreasing with respect to its number of offered products (i.e., given a second-level facility $j\in J$, the cost of offering $t+1$ products is always greater than or equal to the cost of offering $t$ products), and considering that given a second-level facility, it can only be opened at one size by constraints~\eqref{con:nprod2}, it can be stated that: 
\begin{equation}\label{ec:profProp1}
    \sum_{(j,t) \in J_T^1}u_{jt} \leq  \tilde{k},
\end{equation}
where $J_T^1:=\{(j,t)\in J_T:\exists t'\leq t: (j,t')\in J_T'\}$. Since $J_T'\subseteq J_T^1,$ constraint \eqref{ec:profProp1} dominates \eqref{ec:profProp0}. 
Moreover,
\eqref{ec:profProp1} can be again reinforced by including some more $u$-variables in the left-hand side. Let
$$\tilde{c}=\max\left\{B-\sum_{h=1}^{\tilde{k}}\delta_{(h)},\delta_{(\tilde{k})}\right\}.$$ 
Therefore, if we sum to the left-hand side of expression \eqref{ec:profProp1} all the variables $u_{jt}$ with $(j,t)\in J_T$ such that $c_{jt}>\tilde{c}$ and $(j,t)\notin J_T^1$, i.e., $(j,t)\in J_T^2$, the resulting summation will take at most the value $\tilde{k}$. This is valid because we are adding $u$-variables whose cost is greater than the maximum cost associated with $u$-variables included in \eqref{ec:profProp1} and greater than the remaining budget.
From this argument, expression~\eqref{con:DDVVBRef} is derived. 
The above reasoning holds true for any choice of $J_T'\subseteq J_T$, therefore the result is proven. 
\fin

{

In the following, we give a result based on Proposition~\ref{prop:DDVVcovMax} that provides different values of $\tilde{c}$ and $\tilde{k}$ that could result in a reinforcement of constraints~\eqref{con:DDVVBRef}. 

\begin{coro}\label{cor:DDVVcov}
The following family of constraints is valid for the formulation \nameM: 
\begin{equation}\label{con:coro}
 \sum_{(j,t) \in J_T^1 \cup J_T^3} u_{jt} \leq  \bar{k}+\gamma-1,
 \quad  J_T' \subseteq J_T,
 \end{equation}
such that $\gamma\in \mathbb{N}$, $J^1_T$ is defined in Proposition \ref{prop:DDVVcovMax}, $J_T^3:= \left\{(j,t) \in J_T\setminus J_T^1 \: : \: c_{jt} > \bar{c}\right\},$ $\bar{k}$ is the largest index that satisfies the following condition
\begin{equation}\label{con:DDVVB_coro}
    \sum_{h=1}^{\bar{k}} \delta_{(h)}\leq B,
\end{equation}
where $\delta(\cdot)$ is a permutation of costs $c_{jt}$ for 
$(j,t) \in J_T'$
sorted in non-decreasing order, and $\bar{c}=\max\left\{\frac{B-\sum_{h=1}^{\tilde{k}}\delta_{(h)}}{\gamma},\delta_{(\bar{k})}\right\}$.
\end{coro}

\dem 

For the case where $\gamma=1$, the result is the one described in Proposition~\ref{prop:DDVVcovMax}. If $\gamma >1,$ observe that $\gamma-1$ represents the maximum number of $u$-variables with cost at least $\bar{c}$ that can be added to the left-hand side of \eqref{con:DDVVB_coro} without violating this inequality. Note that we are just analyzing the worst-case, which means that it could be the case that there are not $\gamma-1$ $u$-variables with these costs. Hence, in some sense, $\bar{k}+\gamma-1$ is taking the role of $\tilde{k}$ in Proposition~\ref{prop:DDVVcovMax}, where $\delta_{(\bar{k}+1)}=\bar{c},$ ..., $\delta_{(\bar{k}+\gamma-1)}=\bar{c},$ and $\delta_{(\bar{k}+\gamma)}=\bar{c}+\epsilon$, with $\epsilon>0$. 
Moreover,
$$\sum_{h=1}^{\bar{k}+\gamma-1} \delta_{(h)}\leq B< \sum_{h=1}^{\bar{k}+\gamma} \delta_{(h)}$$
 is verified, then $\bar{c}$ is taking the role of $\tilde{c}$ in Proposition~\ref{prop:DDVVcovMax}. Therefore, the result follows.
\fin

Efficient separation procedures for the valid inequalities described in Proposition~\ref{prop:DDVVcovMax} and Corollary~\ref{cor:DDVVcov} are developed. In the following subsection, we detail two separation procedures: one is based on solving auxiliary MILPs and it is designed for the family~\eqref{con:coro} and the other approach chooses the subsets trying to minimize the difference between the right-hand side and the left-hand side of constraints~\eqref{con:DDVVBRef}. 

\subsection{Separation procedures of covering valid inequalities}

As mentioned, the families of constraints \eqref{con:DDVVBRef} and \eqref{con:coro} have an exponential number of constraints. Our goal is to develop a branch and bound and cut procedure to solve SL-MCLP using these constraints. Therefore, for a given fractional solution of the problem, the idea is to identify a constraint of the family \eqref{con:DDVVBRef} (respectively \eqref{con:coro}) that cuts this solution. In this subsection, we propose two different approaches:
\begin{itemize}
    \item The first approach is based on solving a linear integer programming problem that for a given fractional solution of $u$-variables find $J'_T\subseteq J_T$ minimizing the difference between the right and the left-hand side of 
    \eqref{con:coro}.
    \item On the other hand, for a given fractional solution for $u$-variables, the second approach tries to find $J'_T\subseteq J_T$ heuristically minimizing the difference between the right and the left-hand side of \eqref{con:DDVVBRef}.  
    
\end{itemize}  
While the first procedure is computationally more expensive than the second one, the results will be more refined. The advantage of the second procedure is that it is very quick to perform, which is very important since it will have to be repeated at the nodes of the B\&B tree. 

\subsubsection{
Separation procedure of family (\ref{con:coro})} \label{subsec:MIPNew4}
In this subsection, we propose a mathematical formulation for providing valid inequalities of family \eqref{con:coro} that allows us to cut fractional solutions for the $u$-variables in SL-MCFLP. Let $\bar{u}$ be a fractional solution of $u$-variables, 
the idea behind this formulation is to build a constraint of this family minimizing the difference between its right-hand side and the left-hand side. In case this difference is negative, then these constraints will cut the aforementioned solution. Observe that the aim of the program is to identify the elements $(j,t)$ that belong to the set $J_T', J_T^1,$ and $J_T^3$ defined in Proposition~\ref{prop:DDVVcovMax} and Corollary~\ref{cor:DDVVcov}.
Let us consider the following family of binary variables:
\begin{eqnarray*}
\lambda_{jt}& = &
\begin{cases}
1, & \mbox{if $(j,t)\in J_T'$, }\\
0, & \mbox{otherwise,}
\end{cases}\mbox{for $(j,t) \in J_T,$}\\
\mu_{jt}& = &
\begin{cases}
1, & \mbox{if $(j,t)\in J_T^3$,}\\
0, & \mbox{otherwise,} 
\end{cases} \mbox{for $(j,t) \in J_T.$}
\end{eqnarray*}
Therefore, the formulation proposes to identify a valid inequality of the family \eqref{con:coro} is given by:
\begin{eqnarray}
&\min &  \sum_{j\in J} \sum_{t\in T(j)} \lambda_{jt} +(\gamma-1) -\sum_{j\in J} \sum_{t\in T(j)} \lambda_{jt}\left(\sum_{t'= t}^{b_j} \bar{u}_{jt'}\right)-\sum_{j\in J} \sum_{t\in T(j)} \mu_{jt}\bar{u}_{jt}\nonumber\\
\label{sub_c1}
&s.t.&   \sum_{t\in T(j)} \lambda_{jt}\le1, \quad j\in J, \\
\label{sub_c2}
&& \sum_{j\in J} \sum_{t\in T(j)} c_{jt} \lambda_{jt}\le B, \\
\label{sub_c3}
&&  \lambda_{\max} \ge c_{jt}\lambda_{jt}, \quad j\in J, t \in T(j), \\ 
\label{sub_c4}
&& (1-\mu_{jt})B+c_{jt} \ge \lambda_{\max}+\epsilon,  \quad j\in J, t \in T(j), \\
\label{sub_c5}
&& (1-\mu_{jt})B+c_{jt}\gamma \ge B-\sum_{j'\in J}   \sum_{t'\in T(j)} c_{j't'}\lambda_{j't'}+\epsilon,\quad j\in J, t \in T(j),   \\
\label{sub_c6}
&& \mu_{jt} \le 1- \sum_{t' = 1}^{b_j}  \lambda_{jt'}, \quad j\in J, t \in T(j), \\
\label{sub_c7}
&& \gamma \in\mathbb{N}, \\
\label{sub_c8}
&& \lambda_{jt}, \mu_{jt} \in \{0,1\}, \quad j\in J, t \in T(j),\\
&& \lambda_{\max}\geq 0.
\end{eqnarray}

For a given feasible vector value $\bar{u}$ of the $u$-variables, the objective function 
represents the difference between the left-hand side and the right-hand side of a constraint of the family \eqref{con:coro}. Observe that the first two addends of the objective function represent the right-hand side of \eqref{con:coro}. Indeed, the first term represents $\bar{k}$, i.e., a number of $u$-variables whose corresponding costs add up to at most $B$ (this condition is guaranteed by \eqref{sub_c2}). Moreover, constraints \eqref{sub_c1} ensure that at most one size per second-level facility is chosen. 
The third addend of the objective function represents the variables $u_{jt}$ such that $(j,t)\in J_T^1$.  The fourth addend of the objective function represents the variables $u_{jt}$ such that $(j,t)\in J_T^3$. Note that the second and the third addends represent the $u$-variables contained in the set $J_T^1\cup J_T^3$. 
Constraints \eqref{sub_c3} guarantee that $\lambda_{\max}$ takes the maximum cost with value one in the corresponding $\lambda$-variables. By constraints \eqref{sub_c4} and \eqref{sub_c5} with $\epsilon$ a positive constant (close to zero),
 $\mu$-variables represents the $u$-variables with $c$-values bigger than
  $$\max\left\{\frac{B-\sum_{j'\in J}\sum_{t'\in T(j)} c_{j't'}\lambda_{j't'}}{\gamma},\lambda_{\max}\right\}.$$
Constraints \eqref{sub_c6} ensure that each $u$-variable in the valid inequality is chosen either by its corresponding $\lambda$-variable or $\mu$-variable, i.e., belonging to $J_T'$ or $J_T^3$.

Once this model is solved, a valid inequality of family \eqref{con:coro} is obtained and included in the branch and cut procedure.

 
\subsubsection{Separation procedure of family (\ref{con:DDVVBRef})}
\label{subsec:HeuristicNew4} 

In this subsection, we propose a different approach to obtain valid covering inequalities. Using this procedure, we get promising constraints of the family \eqref{con:DDVVBRef}. 
 Therefore, given a fractional solution in a node of the B\&B tree, 
the aim of this section is to select a valid inequality of family 
\eqref{con:DDVVBRef}, i.e., choosing $J_T' \subseteq J_T$ trying to minimize the difference between the right and the left-hand side of \eqref{con:DDVVBRef}. The idea is to find a promising size for each second-level facility to be included in the $J'_T$ set. To do that, we rank the different sizes for each second-level facility based on the cost and the definition of set $J_T^1$. Then, we choose those second-level facilities whose associated $u$-variables are the largest ones, trying to obtain a valid inequality of the family \eqref{con:DDVVBRef} that is not verified at that node. 


\paragraph{First phase}

The objective of this phase is to select for each second-level facility $j\in J$, the most promising size, i.e., $t^\prime(j)\in T(j)$. For that purpose, for each $j\in J$, we obtain: 
$$t^\prime(j)=\text{argmax}_{t\in T(j)}\left\{c_{jt}\sum_{\tilde{t}=t}^{b_j}\bar{u}_{j\tilde{t}}\right\}.$$

\paragraph{Second phase}

The objective of this phase is to select a subset of second-level facilities, i.e., $J^\prime \subseteq J$.

First, we sort the second-level facilities $j \in J$ in non-increasing order according to $\displaystyle\sum_{\tilde{t}=t^\prime (j)}^{b_j}\bar{u}_{j\tilde{t}}$. This means that we prioritize those second-level facilities with the largest $u$-variables values. 

Let $\bar{\delta}(\cdot)$ be a permutation of cost $c_{jt^\prime(j)}$ that sorts the cost in the order described above and $\tilde{k}$ the largest index $k$ that satisfies the following condition:
\begin{equation}\label{con:DDVVAlgHeuPresu}
    \sum_{h=1}^{k} \bar{\delta}_{(h)}\leq B.
\end{equation}

The pair $(j,t'(j))\in J_T$ such that their associated cost ($c_{jt^\prime(j)}$) has been included in constraint \eqref{con:DDVVAlgHeuPresu} are incorporated in $J^\prime_T$. 

\paragraph{Third phase}

Include constraint \eqref{con:DDVVBRef} in the B\&B tree, using $J^\prime_T$ defined above, whenever this constraint cuts the current fractional solution. 

The described procedure has been summarized in Algorithm~\ref{Alg:Heur_pres}.

\begin{algorithm}[tb]\label{Alg:Heur_pres}
\setstretch{1.15}
	\KwData{Let $(\bar{x},\bar{y},\bar{u})$ be a fractional solution of \nameM in a node of the B\&B.}
	\KwResult{The set $J_T'$ required in the valid inequalities~\eqref{con:DDVVBRef}.}
\For{ $j\in J$} {
Obtain $t^\prime(j)=\text{argmax}_{t\in T(j)}\left\{c_{jt}\sum_{\tilde{t}=t}^{b_j}\bar{u}_{j\tilde{t}}\right\}$.

}

Sort $j \in J$ in non-increasing order according to $\displaystyle\sum_{\tilde{t}=t^\prime (j)}^{b_j}\bar{u}_{j\tilde{t}}$.

Let $\bar{\delta}(\cdot)$ be a permutation of cost $c_{jt^\prime(j)}$ that sorts the cost in the order described above and $\tilde{k}$ the largest index $k$ that satisfies the following condition:
\begin{equation} \nonumber
    \sum_{h=1}^{k} \bar{\delta}_{(h)}\leq B.
\end{equation}

Any pair $(j,t'(j))\in J_T$ such that $c_{jt^\prime(j)}$ has been included in the left-hand side of the above inequality is incorporated to $J^\prime_T$. 

Add constraint \eqref{con:DDVVBRef} into the B\&B procedure, using the set $J^\prime_T$ previously defined, whenever this constraint cuts the current fractional solution. 
\caption{Separation procedure of \eqref{con:DDVVBRef}.}
\end{algorithm}
\section{Matheuristic procedure for the \nameM: Kernel Search}\label{sec:Heuristic}

In this section, we propose a matheuristic algorithm to solve formulation \nameM.  This algorithm is inspired by the Adaptive Kernel Search (AKS), originally presented in \cite{AKS}. The Kernel Search based algorithms have been successfully applied in several different problems, as for example, \cite{JanJAnKveZaj21, BalMarRod21,LabMarRod19}, and \cite{LaMaZa22}. The idea of this algorithm is to solve a sequence of mixed-integer problems derived from the original one, where the majority of the binary variables are fixed to zero and only some of them are considered as integers -- this set of integer variables is called Kernel. 


 The heuristic approach that we propose for \nameM has the following characteristics:
\begin{itemize}
\item[i)] An AKS is introduced to deal with $x$-variables.
\item[ii)] The $u$-variables are fixed to zero or considered as binary depending on the $x$-variables in the Kernel.
\end{itemize}
This strategy considers subproblems where the number of binary variables and constraints are reduced with respect to the
original problem. As a consequence, these subproblems have shorter solution times. 

The algorithm has two phases: 

\begin{enumerate}
    \item The first phase determines the initial set of binary $x$-variables; i.e., the initial kernel. Furthermore,  in this step, the remaining $x$-variables are ranked. Note that the most promising variables to take a value of 1 in the optimal solution are considered first. The chosen order has a huge influence on the solution obtained by the heuristic, so we developed three different strategies that are described below. Note that in the original algorithm, several sets of binary variables can be considered, but, in this version, we only apply the Adaptive Kernel Search over $x$-variables, i.e., $u$-variables are adjusted using a different methodology during the procedure.  
    
    \item In the second phase, a sequence of restricted MILP derived from the original problem is solved including constraints that ensure that a progressively better bound on the solution is obtained. The solution and the computational effort required to solve the previously restricted MILP guide the construction of the subsequent Kernels. In each iteration, the Kernel is updated to test other promising variables and to remove useless ones. 
\end{enumerate}

Let $D:=J\times \{m\in M_J(j): j\in J\}$, i.e., the set that represents the indices of  $x$-variables. Given $\mathcal{K}:= \mathcal{K}_j \times \mathcal{K}_m\subseteq D,$ the restricted MILP derived from the original problem will be referred to as $\mbox{\nameM($\mathcal{K}$)},$ where if $(j,m)\in \mathcal{K},$ $x_{jm}$-variable is considered as a binary variable, otherwise it is fixed to zero. For each $j\in \mathcal{K}_j,$ let $\mathcal{T}(j,\mathcal{K})=\{0,\dots, \min\{b_j,|\{m\in \mathcal{K}_m:  (j,m)\in \mathcal{K}\}|\} \}$, i.e, the number of different products that can be served in a second-level facility $j$ according to kernel $\mathcal{K}$.  Furthermore, variables $u_{jt}$, for $j\in \mathcal{K}_j,$ $t\in \mathcal{T}(j,\mathcal{K})$ are considered  as binary variables, all other $u$-variables are fixed to zero.  Moreover, we define $\mathcal{S}_I(i,\mathcal{K})$ as the set of sources $s\in S$ whose product is demanded by client $i$ such that exists $j\in \mathcal{K}_j$ satisfying $d_{ji}\leq r_{m(s)}$ and  $d_{sj}\leq R_{m(s)}$.

Therefore, $\mbox{\nameM($\mathcal{K}$)}$ model is formulated as:
\begin{align}
&\max&&\displaystyle \sum_{i\in I}\sum_{s\in \mathcal{S}_I(i,\mathcal{K})} w_{is}y_{is}\nonumber\\
&\mbox{s.t.}&&\sum_{s\in S_I(i,\mathcal{K})}y_{is}\leq 1,&& i\in I, \\
&&&\sum_{j\in J(i,s): (j,m(s))\in\mathcal{K} }x_{jm(s)}\geq y_{is},&& i\in I, s\in \mathcal{S}_I(i,\mathcal{K}), \\
&&&\sum_{m: (j,m)\in \mathcal{K}}x_{jm}=\sum_{t\in \mathcal{T}(j,\mathcal{K})}t u_{jt},&&j\in \mathcal{K}_j,\\
&&&\sum_{t\in \mathcal{T}(j,\mathcal{K})}u_{jt}= 1,&&j\in \mathcal{K}_j,\\
&&&\sum_{j\in \mathcal{K}_j}\sum_{t\in \mathcal{T}(j,\mathcal{K})}c_{jt} u_{jt}\leq B,\\
&&&0\leq y_{js}\leq 1,&&i\in I,s\in S_I(i,\mathcal{K}), \\
&&&u_{jt}\in \{0,1\},&&j\in \mathcal{K}_j, t\in \mathcal{T}(j,\mathcal{K}), \\
&&&x_{jm}\in \{0,1\},&&(j,m)\in \mathcal{K}. 
\end{align}

In the following subsections, each phase of the procedure is detailed. 

\subsection{First phase}

The aim of this initial phase is to create an initial kernel on $x$-variables and sort the remaining ones according to how likely these variables will take a value of one in the optimal solution. This stage is essential for the success of the algorithm. For this reason, we develop three different strategies to sort the variables. The first strategy is the original strategy of the AKS algorithm. The other two take advantage of the structure of the problem.  

\begin{itemize}
    \item \textbf{Classical (C):} This strategy is the one used in the original AKS algorithm. First, we solve the linear relaxation of \nameM. Let $\breve{x}$ be the optimal values of $x$-variables  and $\breve{r}$ the reduced costs of these variables in the linear relaxation of formulation \nameM. 
 
 The $x$-variables are sorted in non-increasing order with respect to vector $r$, which is computed as:
$$r_{jm}=\begin{cases} 
    \breve{x}_{jm}, &\text{if }\breve{x}_{jm}>0, \\
\breve{r}_{jm},   &\text{otherwise.}
\end{cases}$$
 The initial kernel, $\mathcal{K}_0$, is composed by $(j,m)\in D$ such that $r_{jm}$ is greater than some threshold fixed by the modeler. 
\end{itemize}

Observe that, in this strategy, it could happen that only a subset of $x$-variables associated with a second-level facility $j$ are taken into account in the model. This means that only a strict subset of $T(j)$ can be considered as possible products offered by $j$ in this iteration. The next strategies that we propose differ in this aspect. 

In the remaining strategies,
we sort the potential location $j$ following a certain criterion.
If the location $j$ is selected, all $x_{jm},$ for $m\in T(j)$ are included in the kernel. These two strategies are:

 \begin{itemize}
    
 \item \textbf{Best product/warehouse (BP):} We compute vector $r^{pw}$ as follows:
$$r^{pw}_j=\max_{m\in T(j)}r_{jm}.$$ 
 
We sort the warehouses $j$ in non-increasing order with respect to vector $r^{pw}$. 
 
 \item \textbf{Best warehouse (BW):}  In this method, we calculate vector $r^{w}$ as follows: 
 $$r^{w}_j=\sum_{m\in T(j)}r_{jm}.$$ 
As in the previous method, we sort the warehouses $j$ in non-increasing order with respect to vector $r^{w}$. 
 \end{itemize}
 
 In this two strategies, the initial kernel ($\mathcal{K}_0$) is composed by $(j,m)\in D,$ such that $r_j^{pw}$ (respectively $r_j^{w}$) is greater than some threshold fixed by the modeler.


 Once this initial kernel $\mathcal{K}_0$ is defined, the model $\mbox{\nameM($\mathcal{K}_0$)}$ is solved. Its optimal objective value is a lower bound (LB) of the optimal objective value of \nameM. 
 

\subsection{Second phase}

In this phase, an iterative procedure starts. In each iteration, ($it$), a new set of indices, named $B_{it}\subseteq D$ is added to the kernel. These indices are included in the order determined by vector $r$ (respectively $r^{pw}, r^{w}$), i.e., the most promising $x$-variables to take value one in the optimal solution are considered first. Observe that each index is included at most once. The initial size of $B_{it}$ is a parameter decided by the modeler. In the classical version, the size of this set is the size of $\mathcal{K}_0$, as done in \cite{AKS}. In the warehouse's methods,  we fixed the number of warehouses to include depending on the size of the initial kernel. In the subsequent iterations, the size of $B_{it}$ will be adjusted taking into account the solving time of the previous iteration. 

The model $\mbox{\nameM($\mathcal{K} \cup B_{it}$)}$ is solved adding the following constraints: 
\begin{eqnarray}
&&\sum_{i\in I}\sum_{s\in \mathcal{S}_I(i,\mathcal{K}\cup B_{it})} w_{is}y_{is}\geq \mbox{LB},  \label{betterUB}\\
&&\sum_{(j,m)\in B_{it}}x_{jm}  \geq 1. \label{atleast1}
\end{eqnarray}

The aim of constraint \eqref{betterUB} is to ensure that a better lower bound of the formulation is obtained. Furthermore,  constraint \eqref{atleast1} forces that at least one variable of the set $B_{it}$ is selected, in order to ensure that a different solution is found. If the previous iteration found a feasible solution but the optimality of this solution was not proved because the time limit was reached, constraint \eqref{atleast1} will be replaced  by: 
\begin{equation}
\sum_{(j,m)\in B_{it}  \cup \{(j,m) \in \mathcal{K}: \bar{x}_{jm}=0\}}x_{jm} \geq 1, \label{atleast1-no}
\end{equation}
where $\bar{x}_{jm}$ represents the solution associated with the current lower bound.

A feasible solution of this subproblem  
($\mbox{\nameM($\mathcal{K}\cup B_{it}$)+\eqref{betterUB}+\eqref{atleast1}}$ 
if the previous solution is optimal or $\mbox{\nameM($\mathcal{K}\cup B_{it}$)+\eqref{betterUB}+\eqref{atleast1-no}}$ if the previous solution  is feasible) improves the lower bound. However, this problem may not be feasible, thus, a time limit $t_{limit}$ is imposed. Moreover, we also add a time limit $t_{fea}$ for finding a feasible solution, i.e, the subproblem is stopped if no feasible solution is found within $t_{fea}$. In addition, if the incumbent solution is not improved after a fixed time $t_{inc}$, the subproblem is stopped.

Then, the kernel is updated: the selected variables of the set $B_{it}$ are added to the kernel, denoted as $\mathcal{K}^{+}$, and the variables of the kernel that were not selected in the last $p$ iterations are removed from the kernel, denoted as $\mathcal{K}^-.$ 

In the following iterations, the problem  $\mbox{\nameM($\mathcal{K} \cup B_{it}$)}$ is solved. The solving time of the previous iteration determines the size of  $\mathcal{K} \cup B_{it}$, named $S_{\mathcal{K} \cup B_{it}}$. If the problem is easily solved, i.e., if $t_{\mathcal{K} \cup B_{it}} \leq t_{Easy}$, the number of $x$-variables that are binary in the next step is incremented, i.e., $S_{\mathcal{K} \cup B_{it+1}}=(1+\delta)S_{\mathcal{K} \cup B_{it}}$, where $0\leq \delta \leq 1.$ In the warehouse's methods, we increase the number of second-level facilities to consider in each iteration.

This iterative process is stopped if a criterion defined by the user is fulfilled, for example, doing a fixed number of iterations, establishing a time limit, or fixing a percentage of $x$-variables (warehouses) that should be analyzed.

The proposed heuristic has many parameters: threshold to be included in the initial Kernel, $\delta, t_{limit}, t_{Easy},$ etc. The decision-maker should adapt the values of these parameters to the data and their needs, finding the desired trade-off between time and precision. An overview of the developed Kernel Search procedure is described in Algorithm \ref{a:ks} and Figure \ref{flowchart}.

In conclusion, in this section, we have presented a new algorithm to obtain accurate feasible solutions for our model based on the kernel search algorithm. We have introduced three strategies to sort and choose the most promising $x$-variables. In particular, strategies BP and BW are developed taking advantage of the structure of the problem assuming that, once a location $j$ is selected, all $x_{jm}$ for $m\in T(j)$ are included in the kernel.

\begin{algorithm}[htbp]\label{a:ks}

	\KwData{Instance data of problem \nameM.}
	\KwResult{Heuristic solution of \nameM.}

Solve the linear relaxation of \nameM. \\
Sort $x$-variables using one of the following criteria: classical (C), best product/warehouse (BP), or best warehouse (BW).\\
Constitute  initial $\mathcal{K}=\mathcal{K}_0$ using the first ranked $x$-variables.\\
Solve \nameM($\mathcal{K}$). Let $(\bar{x},\bar{y},\bar{u})$ be the solution of the current LB.

$it=0$

\Do{a criterion is fulfilled (e.g. number of iterations, time limit)}{

Determine a new subset $B_{it}\subseteq D$ following the ranking criterion.

\If{ the solution at $it-1$ is optimal or $it=0$}
{Solve \nameM($\mathcal{K}\cup B_{it}$)+\eqref{betterUB}+\eqref{atleast1}.
}

\If{ the solution at  $it-1$ is feasible but not optimal}
{Solve \nameM($\mathcal{K}\cup B_{it}$)+\eqref{betterUB}+\eqref{atleast1-no}.
}

Let:
\begin{itemize}
\item[]$\mathcal{K}^{+}:=\{k\in B_{it}: \bar{x}_k=1 \}.$ 
\item[]$\mathcal{K}^{-}:=\{k\in \mathcal{K}: k \mbox{ has not been selected in the last $p$ feasible iterations.} \}.$
\end{itemize}
Update $\mathcal{K}=\mathcal{K}\cup\mathcal{K}^{+}\setminus\mathcal{K}^-$. 

Update LB.

Some parameter values of the next iterations are updated depending on how easy or difficult was to solve the present iteration.

 $it=it+1$.
}
\Return $(\bar{x},\bar{y},\bar{u})$
	\caption{Adaptive Kernel Search for \nameM }
\end{algorithm}

\begin{figure}
    \centering
    \includegraphics[width=1\linewidth]{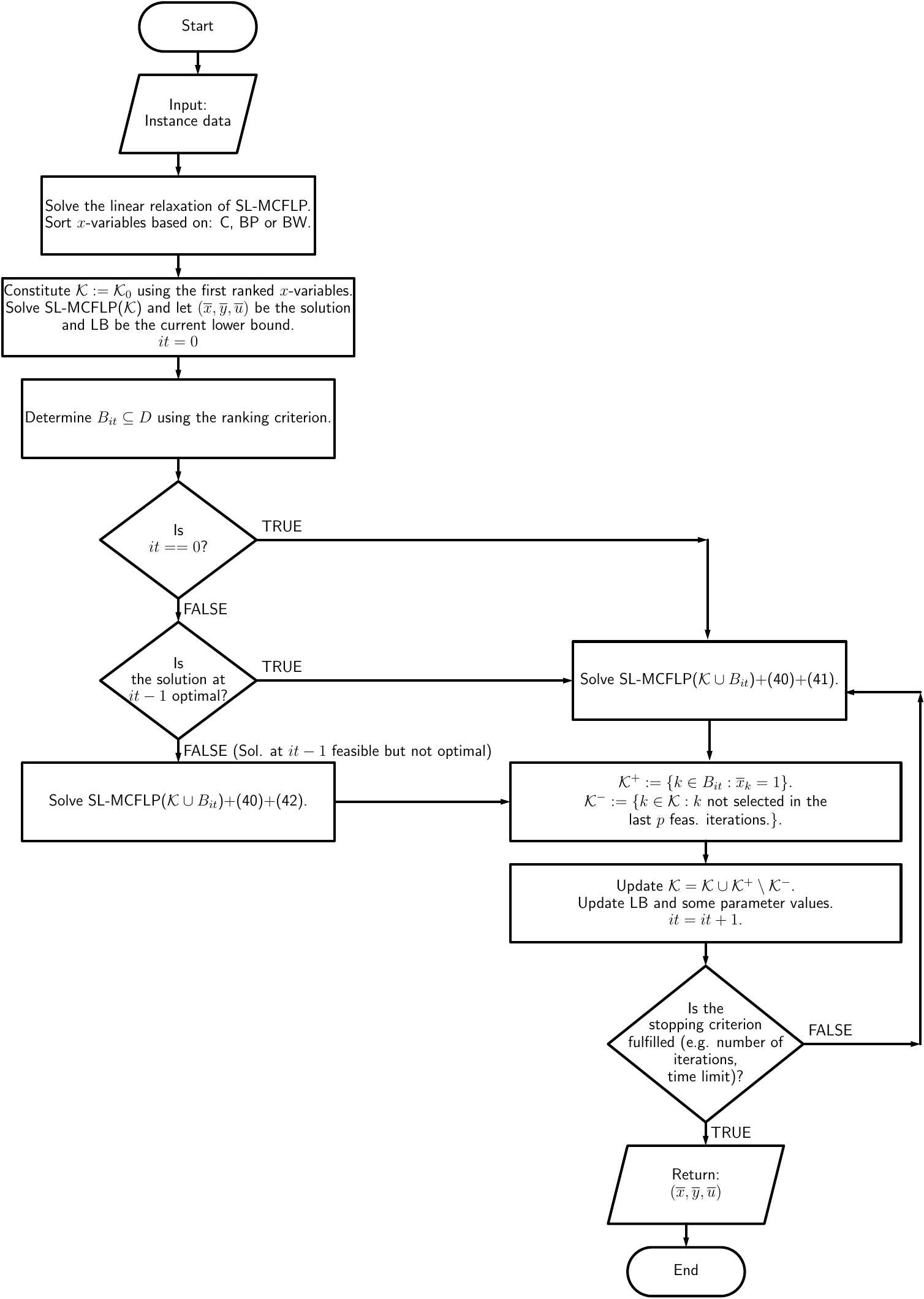}
    \caption{Flowchart of the Adaptive Kernel Search for \nameM.}
    \label{flowchart}
\end{figure}

\section{Computational experience}\label{sec:Computational}

In this section we present two main computational studies. The first one shows the efficiency of the proposed formulation together with the resulting enhancements obtained by including the developed valid inequalities. Secondly, the performance of the proposed variants of the heuristic algorithm is also analyzed.

The experiments were conducted on an Intel(R) Xeon(R) W-2135 CPU 3.70 GHz 32 GB RAM, using CPLEX 20.1.0 in
Concert Technology C++ with a time limit of 7200 seconds. The parameters have been left at their default values where not stated otherwise. 

The instances and the matheuristic codes are available at \url{https://github.com/orel-uca/SLMCFLP}.

\subsection{Instances generation}

For the computational experiments, we considered a different number of first-level facilities: $|S|\in\{30,40\}$. For $|S|=30$, two groups of instances were generated. In the first group, the instances were obtained by creating $|J|=100$ candidate locations for second-level facilities and a different number of clients $|I|\in\{100,250,500\}$. For the second group of instances, we consider $|J|=250$ and $|I|\in \{250,625,1250\}$. Similarly, when $|S|=40$, we studied the cases where $|J|\in\{100,250\}$. For $|S|=40$ and $|J|=100$, the  number of clients were $|I|\in\{100,250,500\}$. Finally, for $|S|=40$ and $|J|=250$, the number of clients varied in $|I|\in\{250,625,1250\}$. 

In each instance, the preference of each client for each first-level facility ($w_{is}$) must be established. For this purpose, we considered two cases: 
\begin{enumerate}
    \item[i)] Direct case ($dir$). For each product $m\in M$, the client prefers the closest first-level facility producing $m$. This configuration tries to describe the behavior of clients who prefer local or zero-kilometer products.
    \item[ii)] Inverse case ($inv$). For each product $m\in M$, the client prefers the farthest first-level facility. Consequently, the client prefers products of exotic origins to the local product.
\end{enumerate}

Observe that, for each instance configuration  the number of products is set to $|M|=\left\lceil \frac{|S|}{4}\right\rceil$.


The locations of first-level facilities, potential second-level facilities, and customers have been randomly generated on a circle of radius 100 centered at the origin. Rounded Euclidean distances are considered among them.
To fix the coverage radii, they have been fixed as the minimum values that guarantee that each first-level facility covers at least one second-level facility, each second-level facility is covered at least by one first-level facility, and covers at least one customer.

Besides, each customer is covered by at least one second-level facility. Given the distribution of the points in the plane, we have observed that this yields instances where some customers have potential access to all first-level facilities, and the case where a customer has potential access to one single first-level facility seldom occurs.
Products are randomly assigned to sources and then some products interchanges are performed to increase the dispersion of the first-level facilities of the same product. 

 For each potential second-level facility, the maximum number of products that can offer is a random number between $2$ and the maximum number of different products whose first-level facilities are in the required radius. Besides, the budget $B$ is set between 10\%-20\% of the costs of offering the first product in all the second-level facilities. The percentage depends on the average number of different products that each second-level facility can offer. It is calculated following the next expression:

$$B=\left\lceil\displaystyle 0.1\cdot \frac{1-0.5^{\left(\frac{\sum_{j\in J}b_j}{|J|}\right)}}{0.5} \cdot \sum_{j\in J}c_{j1}\right\rceil.$$
Observe that $\frac{1-0.5^{\left(\frac{\sum_{j\in J}b_j}{|J|}\right)}}{0.5}$ is a value between 1 and 2. Note that for each configuration, we generated 5 different instances modifying the random seed generator.

\subsection{Performance of the formulation and valid inequalities}

In this subsection, we show the enhancements provided by the inclusion of valid inequalities in formulation \nameM. In preliminary computational experiments, we check the performance of including the families of valid inequalities of polynomial cardinality as constraints in the formulation and in a pool of user cuts. We conclude that the best performance was obtained when the constraints were included directly in the formulation, i.e., as valid inequalities. Additionally, the family of inequalities~\eqref{ec:DDVVcubri} was tested and it was found that it did not provide good computational results.  Regarding the family of constraints    \eqref{ec:DDVVnew4} and \eqref{con:DDVVBRef}, whose size is exponential, we have applied the heuristic separation procedure in the branching tree, described in Algorithm~\ref{strategy1} and Algorithm~\ref{Alg:Heur_pres} 
respectively, using the Generic Callback of CPLEX. In the rest of this paper, unless stated otherwise, whenever the inclusion of these families of constraints is indicated, i.e., \eqref{ec:DDVVnew4} and \eqref{con:DDVVBRef}, it is assumed that these previously mentioned separation procedures have been carried out. Regarding the family of constraints~\eqref{con:coro},  the procedure described in Subsection~\ref{subsec:MIPNew4}  was tested, but the results obtained were worse than those produced by Algorithm~\ref{Alg:Heur_pres} for the family of constraints~\eqref{con:DDVVBRef}. Finally, we would also like to point out that preliminary experiments were carried out with the alternative formulation described in the Appendix~\ref{App:formulation}, \nameMprime. This formulation yielded very similar (although slightly worse) results to the formulation \nameM and this is the reason why they have been omitted in this section.  

In Table \ref{tab:SUM_VI} and Figure \ref{performanceProfilesALL}, a summary of the performance of some combinations of the proposed valid inequalities is detailed. Observe that we have only reported the most promising combinations taking into account that the families of valid inequalities \eqref{ec:DDVVnew1} and \eqref{ec:DDVVnew2} are not considered simultaneously, since the latter is included within the former. The first column of the table reports the formulation used (Formulation). The next column indicates the number of instances solved to optimality within the time limit of two hours (\#Sol.) and the third column, (\%Sol), reports the same information but in percentage. Note that there is a total of 120 instances. Then, the average solution time in seconds is shown in column ``t(s.)". Observe that, for the instances that the optimality was not proven within the time limit, a solution time equal to 7200 s was assigned, i.e.,
the times are underestimated. Then, the average MIP relative gap reported by CPLEX (G\%) is depicted. Finally, it is shown the
linear relaxation gap ($G_{LP}^t$\%), computed as follows:  $$\text{G}_{LP}^{t}\%=\dfrac{\text{ LP}-\text{ BS}^{t}}{\text{ BS}^{t}}\cdot 100,$$ where LP is the optimal solution value of the linear relaxation of the formulation and $\text{BS}^t$ is the best MIP solution value found within the time limit across all formulations. 

For the purpose of a clearer comparison of these formulations, the performance profile graph of the number of solved instances is depicted in Figure~\ref{performanceProfilesALL}. For the sake of clarity, only the formulations with just one family of valid inequalities have been included in this graph. 

\begin{table}[htbp]
  \centering
    \begin{tabular}{lrrrrr}
    Formulation & \multicolumn{1}{l}{\#Sol.} &\multicolumn{1}{l}{\%Sol.}  & \multicolumn{1}{l}{t (s.)} & \multicolumn{1}{l}{$G\%$} & \multicolumn{1}{l}{$\text{G}_{LP}^{t}\%$} \\
    \hline
    \nameM   & 74  &61.67  & 3283.30 & 0.57 & 8.96 \\
    \nameM+\eqref{ec:DDVVpref} & 85 &70.83  & 2424.35 & 0.28 & 8.96 \\
     \nameM+\eqref{ec:DDVVnew4} & 79 &65.83   & 2827.83 & 0.42 & 8.96 \\
    \nameM+\eqref{ec:DDVVnew1} & 90&75.00    & 2180.24 & 0.26 & 3.21 \\
    \nameM+\eqref{ec:DDVVnew2} & 91 &75.83   & 2114.11 & 0.24 & 3.16 \\
    \nameM+\eqref{con:DDVVBRef} & 75 &62.50   & 3259.69 & 0.58 & 8.96 \\
\nameM+\eqref{ec:DDVVpref}+\eqref{ec:DDVVnew1} & 86   &71.67 & 2353.42 & 0.28 & 3.21 \\   
    \nameM+\eqref{ec:DDVVpref}+\eqref{ec:DDVVnew2} & 85 &70.83   & 2392.40 & 0.28 & 3.16 \\
\nameM+\eqref{ec:DDVVpref}+\eqref{ec:DDVVnew1}+\eqref{con:DDVVBRef} & 85 &70.83   & 2347.02 & 0.27 & 3.21 \\\nameM+\eqref{ec:DDVVpref}+\eqref{ec:DDVVnew2}+\eqref{con:DDVVBRef} & 85 &70.83   & 2418.10 & 0.28 & 3.16 \\\nameM+\eqref{ec:DDVVpref}+\eqref{ec:DDVVnew4}+\eqref{ec:DDVVnew1}+\eqref{con:DDVVBRef} & 84  &70.00  & 2531.31 & 0.31 & 3.21 \\
    \nameM+\eqref{ec:DDVVpref}+\eqref{ec:DDVVnew4}+\eqref{ec:DDVVnew2}+\eqref{con:DDVVBRef} & 84 &70.00 & 2562.94 & 0.33 & 3.16 \\
    \hline
    \end{tabular}%
    \caption{Performance of valid inequalities}
  \label{tab:SUM_VI}%
\end{table}%

\begin{figure}[hbt] \centering
	\begin{tikzpicture}[scale=1.1,font=\footnotesize]
		\begin{axis}[axis x line=bottom,  axis y line=left,
			xlabel=$Time(s)$,
			ylabel=\# of solved instances,
			legend style={at={(1.6,0.8)}}]
			\addplot[orange,semithick] plot coordinates {		
(23.066,28)
(23.553,29)
(23.703,30)
(27.705,31)
(29.147,32)
(30.369,33)
(34.955,34)
(40.384,35)
(52.329,36)
(93.838,37)
(96.294,38)
(97.723,39)
(133.492,40)
(144.727,41)
(162.247,42)
(166.381,43)
(215.66,44)
(288.336,45)
(319.313,46)
(437.058,47)
(556.276,48)
(594.957,49)
(657.514,50)
(671.711,51)
(690.065,52)
(694.956,53)
(709.193,54)
(732.461,55)
(777.611,56)
(957.07,57)
(999.756,58)
(1022.41,59)
(1038.64,60)
(1074.99,61)
(1581.56,62)
(1774.07,63)
(2125.17,64)
(2215.35,65)
(2488.75,66)
(3821.29,67)
(4064.7,68)
(4244.19,69)
(4743.51,70)
(4994.84,71)
(5162.06,72)
(5227.44,73)
(6252.9,74)
(7200.23,75)
			};
			\addlegendentry{\tiny\textbf{\nameM}}	
			\addplot[teal,thick, densely dotted] plot coordinates {		
(23.096,40)
(23.253,41)
(24.229,42)
(30.885,43)
(31.452,44)
(42.114,45)
(46.078,46)
(47.619,47)
(50.615,48)
(51.128,49)
(58.29,50)
(59.303,51)
(68.58,52)
(71.24,53)
(71.692,54)
(73.186,55)
(76.191,56)
(77.104,57)
(78.653,58)
(86.041,59)
(91.197,60)
(108.451,61)
(114.975,62)
(125.216,63)
(177.711,64)
(179.051,65)
(245.937,66)
(247.629,67)
(268.03,68)
(418.685,69)
(522.701,70)
(699.379,71)
(767.03,72)
(883.829,73)
(925.804,74)
(1311.04,75)
(1420.24,76)
(1569.16,77)
(1623.2,78)
(1879.51,79)
(3156.96,80)
(3227.2,81)
(3342.57,82)
(3463.55,83)
(3494.08,84)
(7129.98,85)
(7200.21,86)
			};
			\addlegendentry{\tiny\textbf{\nameM + (\ref*{ec:DDVVpref})}}	
   \addplot[red,thick,dash pattern=on 7pt off 1pt] plot coordinates {		

(22.49,21)
(22.831,22)
(23.037,23)
(24.686,24)
(25.778,25)
(28.89,26)
(29.66,27)
(29.739,28)
(30.237,29)
(31.856,30)
(32.953,31)
(35.983,32)
(36.56,33)
(38.539,34)
(41.865,35)
(44.871,36)
(45.269,37)
(46.996,38)
(50.605,39)
(61.915,40)
(67.114,41)
(68.558,42)
(70.939,43)
(77.575,44)
(85.725,45)
(95.14,46)
(96.027,47)
(113.505,48)
(120.294,49)
(123.404,50)
(125.843,51)
(127.979,52)
(145.508,53)
(210.093,54)
(228.449,55)
(244.09,56)
(275.52,57)
(353.702,58)
(359.38,59)
(368.921,60)
(477.802,61)
(516.566,62)
(592.499,63)
(798.633,64)
(801.195,65)
(825.724,66)
(1003.62,67)
(1397.72,68)
(1500.22,69)
(1537.83,70)
(1538.31,71)
(1729.48,72)
(1846.83,73)
(2057.69,74)
(3064.1,75)
(3173.62,76)
(4564.5,77)
(5629.65,78)
(6580.5,79)
(7200.34,80)

			};
			\addlegendentry{\tiny\textbf{\nameM+ (\ref*{ec:DDVVnew4})}}

	\addplot[cyan,thick,dash pattern=on 3pt off 2pt] plot coordinates {		
(24.529,50)
(25.82,51)
(27.292,52)
(28.45,53)
(29.335,54)
(29.639,55)
(32.623,56)
(34.682,57)
(38.38,58)
(44.14,59)
(44.438,60)
(44.914,61)
(47.547,62)
(56.54,63)
(65.585,64)
(71.339,65)
(74.095,66)
(101.357,67)
(115.687,68)
(122.862,69)
(267.634,70)
(377.004,71)
(381.919,72)
(448.336,73)
(469.425,74)
(478.199,75)
(512.128,76)
(525.797,77)
(671.141,78)
(1171.18,79)
(1395.13,80)
(1608.14,81)
(1907.38,82)
(1990.09,83)
(2398.95,84)
(2618.49,85)
(3237.17,86)
(3833.08,87)
(5979.62,88)
(6837.61,89)
(7017.93,90)
(7200.16,91)

			};
			\addlegendentry{\tiny\textbf{\nameM+ (\ref*{ec:DDVVnew1})}}			

	\addplot[magenta,ultra thick, densely dotted] plot coordinates {		
(21.376,50)
(23.577,51)
(23.945,52)
(29.223,53)
(29.515,54)
(33.993,55)
(34.712,56)
(35.067,57)
(37.385,58)
(40.276,59)
(47.058,60)
(48.176,61)
(60.91,62)
(64.215,63)
(65.843,64)
(67.221,65)
(77.444,66)
(83.883,67)
(84.659,68)
(152.558,69)
(154.205,70)
(262.103,71)
(388.512,72)
(430.976,73)
(444.355,74)
(571.547,75)
(576.501,76)
(623.788,77)
(809.11,78)
(942.801,79)
(1046.2,80)
(1174.24,81)
(1189.89,82)
(1296.54,83)
(1634.55,84)
(1643.94,85)
(3463.36,86)
(3562.96,87)
(3610.24,88)
(5930.63,89)
(6563.56,90)
(7105.28,91)
(7200.29,92)

			};
			\addlegendentry{\tiny\textbf{\nameM+ (\ref*{ec:DDVVnew2})}}

   \addplot[violet,thick,densely dashdotted] plot coordinates {		
(20.904,28)
(22.548,29)
(24.735,30)
(25.576,31)
(26.327,32)
(28.759,33)
(29.033,34)
(31.813,35)
(35.226,36)
(47.407,37)
(61.804,38)
(81.857,39)
(114.461,40)
(132.245,41)
(169.36,42)
(181.459,43)
(182.923,44)
(236.173,45)
(290.591,46)
(353.008,47)
(625.078,48)
(635.758,49)
(639.76,50)
(687.271,51)
(793.476,52)
(824.978,53)
(829.213,54)
(835.694,55)
(837.568,56)
(868.765,57)
(1203.06,58)
(1238.73,59)
(1245.13,60)
(1540.68,61)
(1878.47,62)
(1915.61,63)
(2398.81,64)
(2489.58,65)
(2656.05,66)
(2887.67,67)
(3277.87,68)
(3727.07,69)
(3810.84,70)
(4235.81,71)
(5041.18,72)
(5150.05,73)
(5607.38,74)
(6541.69,75)
(7201.7,76)

			};
			\addlegendentry{\tiny\textbf{\nameM+ (\ref*{con:DDVVBRef})}}
		\end{axis}
	\end{tikzpicture}
	\caption{Performance profile graph of \#solved instances using \nameM and the proposed valid inequalities}  \label{performanceProfilesALL}
\end{figure}
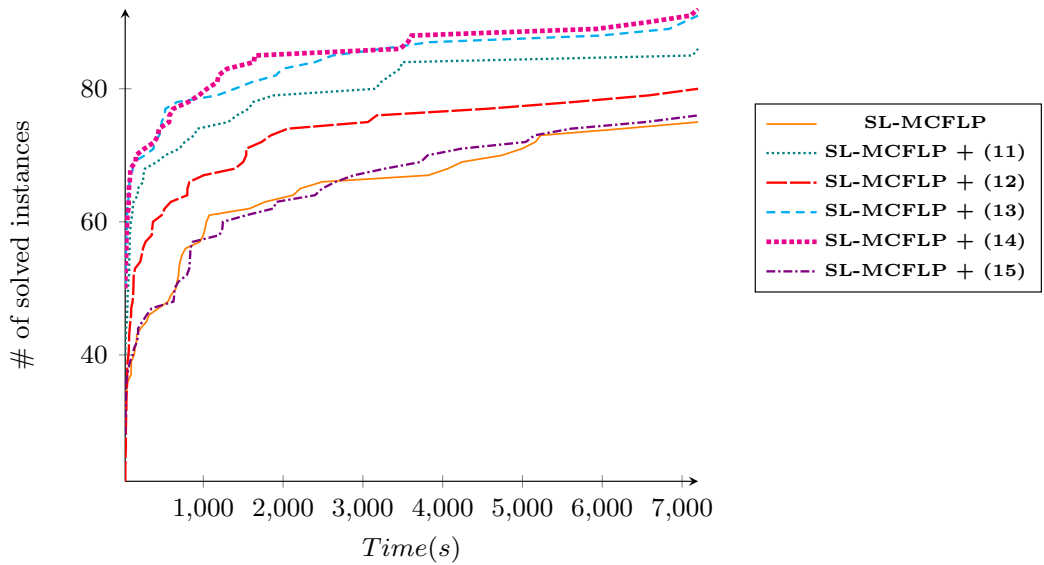

From the above results, it can be inferred that including valid inequalities is very efficient for solving this problem. Not only because they significantly increase the number of instances solved to optimality within the time limit (from 61.67\% to 75.83\% in the best case), but also because they considerably reduce the average MIP relative gap (from 0.57\% to 0.24\% in the best case) and the average gap of the linear relaxation (from 8.96\% to 3.16\% in the best case). Observe that the linear relaxation gap $(G^t_{LP}\%)$ when the family of constraints \eqref{ec:DDVVnew4} (respectively \eqref{con:DDVVBRef}) is included coincides with the one obtained when no constraints are included. This is because these families are included as a separation procedure in the branch and bound tree. 

Next, in Table \ref{tab:Perfnew1new2}, we provide more details on the performance of \nameM and the formulation with the two most promising families of valid inequalities (see Table~\ref{tab:SUM_VI}), i.e., families \eqref{ec:DDVVnew1} and \eqref{ec:DDVVnew2} respectively. This table shows the 
influence on the solving time of the size of the different input sets, as well as the effect of including one of these two families of valid inequalities. The results are the average over five instances generated with the same procedure, varying only the random seed for the generator. The first column indicates the number of first-level facilities ($|S|$). Next, the number of second-level facilities ($|J|$) and clients ($|I|$) are reported, followed by the demand type (D.T.).
The following four columns describe information about \nameM without valid inequalities. The first column reports the number of instances solved to optimality within the time limit. The second one shows the average  time (in seconds) of solving the corresponding five instances. Observe that, if any of these instances is not solved to optimality, 7200 seconds is considered as its solution time to compute this average. Then, the next two columns of this block depict the MIP relative gap reported by CPLEX ($G\%$) and the linear relaxation gap ($\text{G}_{LP}^{t}$\%), respectively. The following blocks of columns report the same information as the other two formulations, i.e., \nameM + \eqref{ec:DDVVnew1}, and \nameM + \eqref{ec:DDVVnew2}.  Finally, the formulation that provided the smallest average solving time is highlighted. If any of the five instances were not solved to optimality, the formulation that solved more instances is shown in bold.

\begin{table}[htbp]
  \centering
  	\resizebox{\linewidth}{!}{
    \begin{tabular}{c|c|c|c|rrrr|rrrr|rrrr}
        \multirow{2}[2]{*}{$|S|$}  &   \multirow{2}[2]{*}{$|J|$}     &   \multirow{2}[2]{*}{$|I|$}    &  \multirow{2}[2]{*}{D.T.}   & \multicolumn{4}{c|}{\nameM}       & \multicolumn{4}{c|}{\nameM + \eqref{ec:DDVVnew1} }    & \multicolumn{4}{c}{\nameM + \eqref{ec:DDVVnew2} } \\
        \cline{5-16}          &       &       &       & \multicolumn{1}{l}{\#Sol.} & \multicolumn{1}{l}{t (s.) } & \multicolumn{1}{l}{G\%} & \multicolumn{1}{l|}{$\text{G}_{LP}^{t}$\%} & \multicolumn{1}{l}{\#Sol.} & \multicolumn{1}{l}{t (s.) } & \multicolumn{1}{l}{G\%} & \multicolumn{1}{l|}{$\text{G}_{LP}^{t}$\% } & \multicolumn{1}{l}{\#Sol.} & \multicolumn{1}{l}{t (s.) } & \multicolumn{1}{l}{G\%} & \multicolumn{1}{l}{$\text{G}_{LP}^{t}$\% } \\
    \hline
    \multirow{12}[5]{*}{30} & \multirow{6}[6]{*}{100} & \multirow{2}[2]{*}{100} &$dir$     & 5     & 3.04  & 0.00 & 9.21 & 5     & 2.00  & 0.00 & 4.39 & \textbf{5} & \textbf{1.81} & \textbf{0.00} & \textbf{4.28} \\
          &       &       & $inv$     & 5     & 10.38 & 0.00 & 5.41 & 5     & 2.19  & 0.00 & 2.28 & \textbf{5} & \textbf{2.09} & \textbf{0.00} & \textbf{2.26} \\
\cline{3-16}          &       & \multirow{2}[2]{*}{250} & {$dir$}     & 5     & 96.22 & 0.00 & 11.66 & 5     & 14.59 & 0.00 & 4.61 & \textbf{5} & \textbf{11.79} & \textbf{0.00} & \textbf{4.52} \\
          &       &       & $inv$     & 5     & 246.55 & 0.00 & 13.15 & 5     & 5.72  & 0.00 & 3.13 & \textbf{5} & \textbf{5.21} & \textbf{0.00} & \textbf{3.12} \\
\cline{3-16}          &       & \multirow{2}[2]{*}{500} & {$dir$}     & 5     & 214.96 & 0.00 & 10.33 & 5     & 18.93 & 0.00 & 3.94 & \textbf{5} & \textbf{15.78} & \textbf{0.00} & \textbf{3.84} \\
          &       &       & $inv$    & 4     & 2107.25 & 0.18 & 7.10 & \textbf{5} & \textbf{91.86} & \textbf{0.00} & \textbf{2.71} & 5     & 209.74 & 0.00 & 2.70 \\
\cline{2-16}          & \multirow{6}[6]{*}{250} & \multirow{2}[2]{*}{250} & $dir$    & 5     & 1127.99 & 0.00 & 8.53 & 5     & 161.61 & 0.00 & 3.20 & \textbf{5} & \textbf{141.70} & \textbf{0.00} & \textbf{3.17} \\
          &       &       & $inv$     & 3     & 5118.44 & 0.17 & 6.47 & 5     & 516.54 & 0.00 & 2.01 & \textbf{5} & \textbf{479.04} & \textbf{0.00} & \textbf{2.01} \\
\cline{3-16}          &       & \multirow{2}[2]{*}{625} & $dir$     & 1     & 6263.31 & 1.05 & 10.30 & 2     & 4824.65 & 0.54 & 3.50 & \textbf{3} & \textbf{4655.42} & \textbf{0.40} & \textbf{3.41} \\
          &       &       & $inv$    & 0     & 7208.40 & 1.25 & 6.40 & 2     & 5326.19 & 0.33 & 2.50 & \textbf{2} & \textbf{5196.10} & \textbf{0.31} & \textbf{2.50} \\
\cline{3-16}          &       & \multirow{2}[2]{*}{1250} & $dir$& 0     & 7203.42 & 1.92 & 11.19 & 2     & 6372.65 & 0.80 & 3.96 & \textbf{3} & \textbf{6102.46} & \textbf{0.85} & \textbf{3.85} \\
          &       &       & $inv$     & 0     & 7204.83 & 2.31 & 6.10 & 0     & 7205.20 & 1.16 & 2.34 & 0     & 7204.51 & 1.47 & 2.33 \\
    \hline
    \multirow{12}[5]{*}{40} & \multirow{6}[6]{*}{100} & \multirow{2}[2]{*}{100} & {$dir$}     & 5     & 41.82 & 0.00 & 9.61 & \textbf{5} & \textbf{4.72} & \textbf{0.00} & \textbf{3.48} & 5     & 5.09  & 0.00 & 3.34 \\
          &       &       & $inv$     & 5     & 12.85 & 0.00 & 9.75 & 5     & 2.31  & 0.00 & 2.81 & \textbf{5} & \textbf{2.15} & \textbf{0.00} & \textbf{2.80} \\
\cline{3-16}          &       & \multirow{2}[2]{*}{250} & $dir$    & 5     & 860.35 & 0.00 & 14.90 & 5     & 30.96 & 0.00 & 4.79 & \textbf{5} & \textbf{13.43} & \textbf{0.00} & \textbf{4.73} \\
          &       &       & $inv$     & 5     & 286.59 & 0.00 & 9.11 & 5     & 5.24  & 0.00 & 2.54 & \textbf{5} & \textbf{5.22} & \textbf{0.00} & \textbf{2.53} \\
\cline{3-16}          &       & \multirow{2}[2]{*}{500} & $dir$    & 4     & 2946.48 & 0.23 & 11.32 & 5     & 356.35 & 0.00 & 4.76 & \textbf{5} & \textbf{148.80} & \textbf{0.00} & \textbf{4.62} \\
          &       &       & $inv$     & 3     & 3216.32 & 0.22 & 7.90 & \textbf{5} & \textbf{20.24} & \textbf{0.00} & \textbf{2.80} & 5     & 39.34 & 0.00 & 2.79 \\
\cline{2-16}          & \multirow{6}[6]{*}{250} & \multirow{2}[2]{*}{250} & $dir$    & 4     & 3021.56 & 0.14 & 7.80 & \textbf{4} & \textbf{1781.12} & \textbf{0.06} & \textbf{3.35} & 4     & 1924.22 & 0.04 & 3.27 \\
          &       &       & $inv$     & 3     & 4101.85 & 0.23 & 8.34 & \textbf{5} & \textbf{926.44} & \textbf{0.00} & \textbf{1.84} & 4     & 1589.87 & 0.03 & 1.84 \\
\cline{3-16}          &       & \multirow{2}[2]{*}{625} & $dir$    & 2     & 5888.68 & 0.37 & 9.04 & 4     & 3115.78 & 0.14 & 4.22 & \textbf{4} & \textbf{2579.81} & \textbf{0.16} & \textbf{4.17} \\
          &       &       & $inv$     & 0     & 7208.24 & 1.16 & 4.53 & 0     & 7206.93 & 0.57 & 2.05 & \textbf{0} & \textbf{7206.33} & \textbf{0.54} & \textbf{2.04} \\
\cline{3-16}          &       & \multirow{2}[2]{*}{1250} & $dir$    & 0     & 7204.20 & 2.11 & 9.96 & 1     & 7129.38 & 1.16 & 3.38 & \textbf{1} & \textbf{5996.24} & \textbf{0.80} & \textbf{3.26} \\
          &       &       & $inv$     & 0     & 7205.55 & 2.34 & 6.84 & 0     & 7204.09 & 1.30 & 2.53 & \textbf{0}     & \textbf{7202.44} & \textbf{1.16} & \textbf{2.52} \\
    \hline
    \end{tabular}%
    }
    \caption{Performance of formulation \nameM, \nameM + \eqref{ec:DDVVnew1}, and \nameM + \eqref{ec:DDVVnew2}.}
  \label{tab:Perfnew1new2}%
\end{table}%

Table \ref{tab:Perfnew1new2} shows that \nameM is the formulation with the most limited results in terms of performance. Observe that, for $|S|=\{30,40\}$, $|J|=250$ and $|I|=1250$ none of the instances are solved in less than two hours since these are the most difficult instance sizes to solve.
 Clearly, formulations  \nameM~+~\eqref{ec:DDVVnew1} and \nameM~+~\eqref{ec:DDVVnew2} outperform \nameM, as previously noted. It can also be appreciated that in the majority of cases, 79.17\%, the family of valid inequalities \eqref{ec:DDVVnew2} is more efficient. However, depending on the sizes of the input sets, there are instances in which it is more efficient to include the family of valid inequalities \eqref{ec:DDVVnew1} rather than the family of valid inequalities \eqref{ec:DDVVnew2}. Note that, in terms of gap percentage after limit time, formulations \nameM~+~\eqref{ec:DDVVnew1} and \nameM~+~\eqref{ec:DDVVnew2} report better results than formulation \nameM.
The details of the remaining formulations analyzed in Table~\ref{tab:SUM_VI} have been omitted for lack of space, but are available on request. 

\subsection{Performance of the heuristic algorithms}
In this section, we validate the efficiency of the different variants of Algorithm \ref{a:ks}, i.e, we compare the solution given by the previous formulations to the ones obtained using the different variants of Algorithm \ref{a:ks}. In some preliminary computational experiments, we checked the efficiency of the Algorithm \ref{a:ks} without and with valid inequalities, achieving considerably better results when valid inequalities were included. From this study, it has been decided to use the family of valid inequalities~\eqref{ec:DDVVnew2}.

We used CPLEX with default values, except for the following parameters:
BndStrenInd = 1, MIPEmphasis = HiddenFeas, FPHeur = 1, as suggested in \cite{AKS}. After a preliminary computational experiment, parameter $\delta$ was set to 0.35, tEasy = 10, tfea = 120, and tinc = 100. Furthermore, in the restricted problems, we limited the computational time (tlimit) to 200 seconds. 

The summary of the results is presented in Table~\ref{tab:summaryheuristic}. The first column indicates the
name of the variant and the second column the average time in seconds. 
Finally, the last column reports the  average percentage difference between
the best solution found by the heuristic algorithm ($\text{BS}^h$) and the best solution found by the exact solution method ($\text{BS}^t$). This difference is
computed as follows: $$\text{G}_{BS}\%=\dfrac{\text{ BS}^t-\text{ BS}^h}{\text{ BS}^t}\cdot 100.$$ Recall that $\text{BS}^t$ is the best MIP solution value found within the time limit across all formulations. 

\begin{table}[htbp]
  \centering
  
    \begin{tabular}{l|rr}
    Variant & \multicolumn{1}{l}{t (s.)} 
    & \multicolumn{1}{l}{$\text{G}_{BS}$\%} \\
    \hline
    C     & 64.72 
    & 0.95 \\
    BP    & 128.72 
    & 0.21 \\
    BW    & 204.56 
    & 0.19 \\
    \hline
    \end{tabular}%
    
     \caption{Summary of the performance of the variants of the heuristic Algorithm. }
  \label{tab:summaryheuristic}%
\end{table}%

 In Table 5, we can appreciate that the average gap percentage with respect the solution of the exact method is smaller than a $1\%$ for all variants of the matheuristic. Besides, the BW variant seems to provide the best quality solutions obtaining an average gap of $0.19\%$. Note that, although BW approach obtains the best average gap, this variant is the most demanding in terms of time. The variant providing the best time results is the classical one (C). 

Furthermore, a raincloud plot of the $\text{G}_{BS}$ and the times of the different variants of the heuristic are depicted in Figures \ref{fig:GAPHueristic} and \ref{fig:timeHueristic}, respectively. The raincloud plot reports a jitter plot, a boxplot (the horizontal bar represents the median), and  an illustration of the data distribution,  for each variant of the heuristic. This graph provides a more comprehensive analysis of the differences in performance between the proposed variants. 

\begin{figure}[htbp]
\centering
\includegraphics[width=\textwidth]{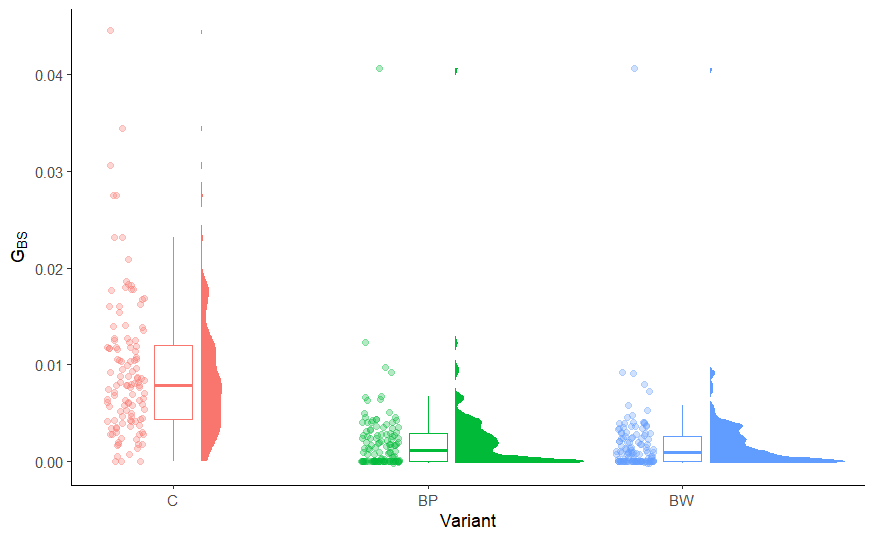}
\caption{Raincloud plot of the GAP $\text{G}_{BS}$ of the different variants of the heuristic Algorithm.}
 \label{fig:GAPHueristic}%
\end{figure}

\begin{figure}[htbp]
\centering
\includegraphics[width=\textwidth]{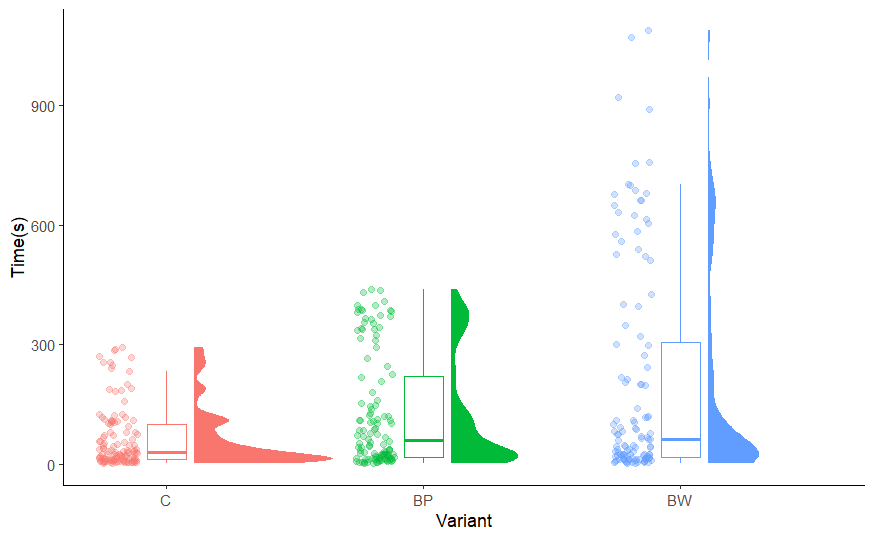}
\caption{Raincloud plot of the time (s) of the different variants of the heuristic Algorithm.}
\label{fig:timeHueristic}
\end{figure}

Considering the reported results, it can be stated that the new variants of the algorithm presented, although 
 with a slightly higher time consumption, provide  better solutions than the classical variant. This can be clearly seen in Figure~\ref{fig:GAPHueristic}, where the GAP ($\text{G}_{BS}$) achieved by the different variants in all instances are plotted. On the other hand, definitely, the heuristic algorithm is efficient, as it provides almost the same solution as the exact algorithm in much less time. Hence, it can be observed that, on average, the best-performing combination of valid inequalities with the formulation (\nameM+\eqref{ec:DDVVnew2}) took 2114.11 seconds, while the slowest version of the heuristic requires a little more than 200 seconds. In fact, in some instances (two for the BP and three for the BW variant), the heuristic algorithm has provided a slightly better solution than the exact algorithm in two hours. Note that this can happen in those instances where the best solution found by all the versions of the formulation in two hours could not be proven to be optimal, i.e., 28 instances.  

Next, Table~\ref{tab:detalleHeu} provides a detailed comparison of the proposed variants taking into account the size of the input sets. The structure of the table is similar to the ones described previously. The variant that provided the smallest best solution GAP ($\text{G}_{BS}\%$) is shown in bold. In the event of a tie, the one with the shortest computational time is highlighted. 
\begin{table}[htbp]
  \centering
  \resizebox{\linewidth}{!}{
    \begin{tabular}{c|c|c|c|rr|rr|rr}
    \multirow{2}[3]{*}{$|S|$} & \multirow{2}[3]{*}{$|J|$} & \multirow{2}[3]{*}{$|I|$} & \multirow{2}[3]{*}{D.T.} & \multicolumn{2}{c|}{C} & \multicolumn{2}{c|}{BP} & \multicolumn{2}{c}{BW} \\
\cline{5-10}          &       &       &       & \multicolumn{1}{l}{t(s.) }  & \multicolumn{1}{l|}{$\text{G}_{BS}\%$ } & \multicolumn{1}{l}{t(s.) }  & \multicolumn{1}{l|}{$\text{G}_{BS}\%$} & \multicolumn{1}{l}{t(s.) }  & \multicolumn{1}{l}{$\text{G}_{BS}\%$ } \\
    \hline
    \multirow{12}[12]{*}{30} & \multirow{6}[6]{*}{100} & \multirow{2}[2]{*}{100} & {$dir$}    & 5.15   & 0.51 & 7.61   & 0.23 & \textbf{8.75} &\textbf{0.15} \\
          &       &       & {$inv$}     & 5.27  & 0.64 & \textbf{6.52} &  \textbf{0.02} & 7.31  & 0.02 \\
\cline{3-10}          &       & \multirow{2}[2]{*}{250} & {$dir$}     & 14.83  & 1.68 & 22.19  & 0.05 & \textbf{22.55}  & \textbf{0.00} \\
          &       &       & {$inv$}     & 11.26  & 0.84 & \textbf{16.91} & \textbf{0.16} & 19.95 &  0.18 \\
\cline{3-10}          &       & \multirow{2}[2]{*}{500} & {$dir$}     & 18.85& 0.95 & \textbf{40.39} & \textbf{0.20} & 45.57  & 0.25 \\
          &       &       & {$inv$}     & 23.44  & 0.61 & 47.83  & 0.25 & \textbf{54.93} &  \textbf{0.05} \\
\cline{2-10}          & \multirow{6}[6]{*}{250} & \multirow{2}[2]{*}{250} & {$dir$}     & 27.35  & 0.60 & 40.75  & 0.15 & \textbf{41.49}  & \textbf{0.12} \\
          &       &       & { $inv$}     & 40.27  & 0.55 & \textbf{102.98}  & \textbf{0.14} & 69.49 &  0.15 \\
\cline{3-10}          &       & \multirow{2}[2]{*}{625} & {$dir$}     & 82.21  & 1.36 & 239.90  & 0.36 & \textbf{403.65}  & \textbf{0.28} \\
          &       &       & { $inv$}     & 109.14  & 0.74 & 185.70 & 0.33 & \textbf{303.29} &  \textbf{0.20} \\
\cline{3-10}          &       & \multirow{2}[2]{*}{1250} & {$dir$}     & 138.31  & 2.01 & \textbf{329.77} & \textbf{0.27} & 464.70  & 0.27 \\
          &       &       & { $inv$}     & 234.29  & 1.24 & 351.12  & 0.56 & \textbf{697.21} &  \textbf{0.51} \\
    \hline
    \multirow{12}[12]{*}{40} & \multirow{6}[6]{*}{100} & \multirow{2}[2]{*}{100} & {$dir$}    & 9.29  & 0.70 & \textbf{16.03} & \textbf{0.00} & 13.33 &  0.03 \\
          &       &       & { $inv$}     & 5.19   & 0.44 & 9.33  & 0.08 & \textbf{9.06} & \textbf{0.00} \\
\cline{3-10}          &       & \multirow{2}[2]{*}{250} & {$dir$}    & 12.22  & 0.62 & \textbf{40.18} &  \textbf{0.02} & 29.39  & 0.20 \\
          &       &       &{ $inv$}     & 13.46  & 1.09 & 15.28  & 0.13 & \textbf{14.48} &  \textbf{0.08} \\
\cline{3-10}          &       & \multirow{2}[2]{*}{500} & {$dir$}     & 48.51  & 1.23 & 108.04  & 0.11 & \textbf{101.74}  & \textbf{0.07} \\
          &       &       & { $inv$}     & 36.38  & 0.63 & 51.28  & 0.06 & \textbf{49.22} & \textbf{0.05} \\
\cline{2-10}          & \multirow{6}[6]{*}{250} & \multirow{2}[2]{*}{250} & {$dir$}     & 35.81  & 0.49 & 172.37 & 0.15 & \textbf{229.31}  & \textbf{0.14} \\
          &       &       & { $inv$}     & 53.89  & 0.39 & \textbf{80.84}  & \textbf{0.05} & 150.90 & 0.08 \\
\cline{3-10}          &       & \multirow{2}[2]{*}{625} & {$dir$}     & 84.30  & 0.86 & \textbf{213.34} &  \textbf{0.26} & 315.95  & 0.28 \\
          &       &       & { $inv$}    & 119.58 & 0.71 & 274.54 & 0.10 & \textbf{502.31} &  \textbf{0.09} \\
\cline{3-10}          &       & \multirow{2}[2]{*}{1250} & {$dir$}     & 178.61 & 2.56 & \textbf{353.72} &  \textbf{1.04} & 631.96 & 1.05 \\
          &       &       & { $inv$}    & 245.69 & 1.30 & 362.58  & 0.40 & \textbf{722.88} & \textbf{0.32} \\
    \hline
    \end{tabular}%
    }
      \caption{Performance of the variants of the heuristic algorithm.}
  \label{tab:detalleHeu}%
\end{table}%

Table \ref{tab:detalleHeu} illustrates that in all cases, the solutions provided by variants BP and BW of the heuristic algorithm are better than the classical variant. Furthermore, it is worth emphasizing that the heuristic algorithm finds the optimal solution or nearly in many cases, given that the GAP with the best solution $\text{G}_{BS}\%$ is extremely small, and in most cases (92/120) this best solution has been proven to be the optimal one. In addition, it is possible to compare the heuristic resolution times of the different instances with the exact resolution time (see Table~\ref{tab:Perfnew1new2}). As an example, we can examine the case of the parameters $|S|=30, |J|=250, |I|=1250$ D.T. equal to d, where the exact algorithm requires 6102.46 seconds on average, (note that 2 of the 5 instances are not solved to optimality, so the times are underestimated) whereas the heuristic algorithm uses only 329.77 seconds and obtains similar solutions ($G_{BS}=0.27\%$ on average). Consequently, it has been shown that the proposed heuristic algorithm is efficient for obtaining good solutions for the problem, as well as, that the new proposed variants perform better than the classical variant.


 \section{Concluding remarks} \label{sec:conclusion}



In this work, the multi-product maximal covering second-level facility location problem (\nameM) is introduced. This model aims to locate a set of second-level facilities, to decide their sizes and the products to be offered, in such a way that the overall clients' satisfaction with respect to their coverage is maximized. To the best of our knowledge, this is the first model for the location of second-level facilities 
which includes clients' preferences in a multi-product and single-flow covering context. In order to satisfy a customer's demand in this model, there must be a double coverage, the customer must be covered by a second-level facility, and this, in turn, by a first-level facility.  

To deal with this model, a MILP is introduced and reinforced with different families of valid inequalities. Some of them are exponential. For this reason, separation algorithms and branch and cut approaches are developed. In the computational tests, we can observe that these reinforcements enhance the resolution of the formulation, increasing the number of solved instances in two hours. 
In particular, constraints \eqref{ec:DDVVnew1} and \eqref{ec:DDVVnew2} are the ones that provide the best results. Besides, three variants of a matheuristic based on Kernel Search are presented. The last two variants, (BP) and (BW), exploit the structure of the model and provide the best results in terms of solution quality.
}

 The characteristics of this model make it useful for e-commerce issues where demand must be delivered in a short period of time. The application of these models will lead to savings in costs, reduction in delivery times, and improved customer satisfaction. We believe this paper provides an insightful starting point for future models. However, from a practical point of view, there may be specific applications in which  the first-level facilities have production limitation. Moreover,  in some applications, the demand is uncertain.
In addition, sometimes it is possible to reduce transport times or distances between two locations for a fixed fee. Therefore, starting from this study, there are several open lines of future research, like adding capacities to the facilities, including uncertainty in demand, 
or enabling 
distance upgrades between two nodes for a fee.

\section*{Acknowledgements}
Research partially supported by: the Spanish Ministry of Science and Innovation through project RED2022-134149-T; Agencia Estatal de Investigación, Spain and ERDF through project PID2020-114594GBC22; MCIN/AEI/10.13039/501100011033 and the European Union ``NextGenerationEU''/PRTR through project TED2021-130875B-I00. Luisa Martínez also acknowledges project AT21 00032 from Junta de Andalucía, Spain and European Union, and project OPTICODS (CIGE/2021/161) from Generalitat Valenciana, Spain.

We would like to thank Maria Albareda Sambola for her collaboration in the early stages of this research.


\begin{thebibliography}{31}
\providecommand{\natexlab}[1]{#1}
\providecommand{\url}[1]{\texttt{#1}}
\expandafter\ifx\csname urlstyle\endcsname\relax
  \providecommand{\doi}[1]{doi: #1}\else
  \providecommand{\doi}{doi: \begingroup \urlstyle{rm}\Url}\fi

\bibitem[Afshari and Peng(2014)]{afshari2014challenges}
H.~Afshari and Q.~Peng.
\newblock Challenges and solutions for location of healthcare facilities.
\newblock \emph{Industrial Engineering and Management}, 3\penalty0
  (2):\penalty0 1--12, 2014.

\bibitem[Baldomero-Naranjo et~al.(2021{\natexlab{a}})Baldomero-Naranjo,
  Kalcsics, and Rodr\'iguez-Ch\'ia]{BalKalRod20}
M.~Baldomero-Naranjo, J.~Kalcsics, and A.~M. Rodr\'iguez-Ch\'ia.
\newblock Minmax regret maximal covering location problems with edge demands.
\newblock \emph{Computers \& Operations Research}, 130:\penalty0 105181,
  2021{\natexlab{a}}.

\bibitem[Baldomero-Naranjo et~al.(2021{\natexlab{b}})Baldomero-Naranjo,
  Martínez-Merino, and Rodríguez-Chía]{BalMarRod21}
M.~Baldomero-Naranjo, L.~I. Martínez-Merino, and A.~M. Rodríguez-Chía.
\newblock A robust {SVM}-based approach with feature selection and outliers
  detection for classification problems.
\newblock \emph{Expert Systems with Applications}, 178:\penalty0 115017,
  2021{\natexlab{b}}.

\bibitem[Baldomero-Naranjo et~al.(2022)Baldomero-Naranjo, Kalcsics, Marín, and
  Rodr\'iguez-Ch\'ia]{BalKalMarRod22}
M.~Baldomero-Naranjo, J.~Kalcsics, A.~Marín, and A.~M. Rodr\'iguez-Ch\'ia.
\newblock Upgrading edges in the maximal covering location problem.
\newblock \emph{European Journal of Operational Research}, 303\penalty0
  (1):\penalty0 14--36, 2022.

\bibitem[Blanco et~al.(2023)Blanco, Gázquez, and Saldanha-da
  Gama]{BlaGazSad23}
V.~Blanco, R.~Gázquez, and F.~Saldanha-da Gama.
\newblock Multi-type maximal covering location problems: Hybridizing discrete
  and continuous problems.
\newblock \emph{European Journal of Operational Research}, 307\penalty0
  (3):\penalty0 1040--1054, 2023.

\bibitem[Casas-Ramírez et~al.(2020)Casas-Ramírez, Camacho-Vallejo, Díaz, and
  Luna]{Casas2020}
M.~Casas-Ramírez, J.~Camacho-Vallejo, J.~A. Díaz, and D.~E. Luna.
\newblock A bi-level maximal covering location problem.
\newblock \emph{Operational Research}, 20\penalty0 (2):\penalty0 827--855,
  2020.

\bibitem[Chauhan et~al.(2019)Chauhan, Unnikrishnan, and
  Figliozzi]{CHAUHAN20191}
D.~Chauhan, A.~Unnikrishnan, and M.~Figliozzi.
\newblock Maximum coverage capacitated facility location problem with range
  constrained drones.
\newblock \emph{Transportation Research Part C: Emerging Technologies},
  99:\penalty0 1--18, 2019.

\bibitem[Church and Murray(2018)]{ChurchMurray2018}
R.~Church and A.~Murray.
\newblock Extended forms of coverage.
\newblock In \emph{Location Covering Models}, pages 49--79. Advances in Spatial
  Science, Springer International Publishing, 2018.

\bibitem[Church and ReVelle(1974)]{ChuRev74}
R.~Church and C.~ReVelle.
\newblock The maximal covering location problem.
\newblock \emph{Papers of the Regional Science Association}, 3:\penalty0
  101--118, 1974.

\bibitem[Contreras and Ortiz-Astorquiza(2019)]{Contreras2019}
I.~Contreras and C.~Ortiz-Astorquiza.
\newblock Hierarchical facility location problems.
\newblock In G.~Laporte, S.~Nickel, and F.~Saldanha~da Gama, editors,
  \emph{Location Science}, pages 365--389. Springer International Publishing,
  2019.

\bibitem[Domínguez et~al.(2021)Domínguez, Labbé, and Marín]{DOLAMA2021}
C.~Domínguez, M.~Labbé, and A.~Marín.
\newblock The rank pricing problem with ties.
\newblock \emph{European Journal of Operational Research}, 294\penalty0
  (2):\penalty0 492--506, 2021.

\bibitem[Farahani et~al.(2014)Farahani, Hekmatfar, Fahimnia, and
  Kazemzadeh]{reviewHierarchical14}
R.~Z. Farahani, M.~Hekmatfar, B.~Fahimnia, and N.~Kazemzadeh.
\newblock Hierarchical facility location problem: Models, classifications,
  techniques, and applications.
\newblock \emph{Computers \& Industrial Engineering}, 68:\penalty0 104--117,
  2014.

\bibitem[Fröhlich et~al.(2020)Fröhlich, Maier, and Hamacher]{FroMaiHam20}
N.~Fröhlich, A.~Maier, and H.~W. Hamacher.
\newblock Covering edges in networks.
\newblock \emph{Networks}, 75\penalty0 (3):\penalty0 278--290, 2020.

\bibitem[Garc{\'i}a and Mar{\'i}n(2019)]{GarMar19}
S.~Garc{\'i}a and A.~Mar{\'i}n.
\newblock Covering location problems.
\newblock In G.~Laporte, S.~Nickel, and F.~Saldanha~da Gama, editors,
  \emph{Location Science}, pages 99--119. Springer International Publishing,
  2019.

\bibitem[Guastaroba et~al.(2017)Guastaroba, Savelsbergh, and Speranza]{AKS}
G.~Guastaroba, M.~Savelsbergh, and M.~Speranza.
\newblock Adaptive kernel search: A heuristic for solving mixed integer linear
  programs.
\newblock \emph{European Journal of Operational Research}, 263\penalty0
  (3):\penalty0 789--804, 2017.

\bibitem[Huang et~al.(2015)Huang, Wang, Batta, and Nagi]{HUANG2015}
S.~Huang, Q.~Wang, R.~Batta, and R.~Nagi.
\newblock An integrated model for site selection and space determination of
  warehouses.
\newblock \emph{Computers \& Operations Research}, 62:\penalty0 169--176, 2015.

\bibitem[Jáno{\v{s}}íková et~al.(2021)Jáno{\v{s}}íková, Jankovič, Kvet,
  and Zajacová]{JanJAnKveZaj21}
L.~Jáno{\v{s}}íková, P.~Jankovič, M.~Kvet, and F.~Zajacová.
\newblock Coverage versus response time objectives in ambulance location.
\newblock \emph{International Journal of Health Geographics}, 20\penalty0 (32),
  2021.

\bibitem[Labbé et~al.(2019)Labbé, Martínez-Merino, and
  Rodríguez-Chía]{LabMarRod19}
M.~Labbé, L.~I. Martínez-Merino, and A.~M. Rodríguez-Chía.
\newblock Mixed integer linear programming for feature selection in support
  vector machine.
\newblock \emph{Discrete Applied Mathematics}, 261:\penalty0 276--304, 2019.

\bibitem[Lamanna et~al.(2022)Lamanna, Mansini, and Zanotti]{LaMaZa22}
L.~Lamanna, R.~Mansini, and R.~Zanotti.
\newblock A two-phase kernel search variant for the multidimensional
  multiple-choice knapsack problem.
\newblock \emph{European Journal of Operational Research}, 297\penalty0
  (1):\penalty0 53--65, 2022.

\bibitem[Landete and Laporte(2019)]{Landete2019}
M.~Landete and G.~Laporte.
\newblock Facility location problems with user cooperation.
\newblock \emph{TOP}, 27\penalty0 (1):\penalty0 125--145, 2019.

\bibitem[Li et~al.(2018)Li, Ramshani, and Huang]{li2018cooperative}
X.~Li, M.~Ramshani, and Y.~Huang.
\newblock Cooperative maximal covering models for humanitarian relief chain
  management.
\newblock \emph{Computers \& Industrial Engineering}, 119:\penalty0 301--308,
  2018.

\bibitem[Mar\'in et~al.(2018)Mar\'in, Mart\'inez-Merino, Rodr\'iguez-Ch\'ia,
  and da~Gama]{MarMarRodSal18}
A.~Mar\'in, L.~I. Mart\'inez-Merino, A.~M. Rodr\'iguez-Ch\'ia, and F.~S.
  da~Gama.
\newblock Multi-period stochastic covering location problems: Modeling
  framework and solution approach.
\newblock \emph{European Journal of Operational Research}, 268\penalty0
  (2):\penalty0 432 -- 449, 2018.

\bibitem[Millstein et~al.(2022)Millstein, Bilir, and Campbell]{MILLSTEIN2022}
M.~A. Millstein, C.~Bilir, and J.~F. Campbell.
\newblock The effect of optimizing warehouse locations on omnichannel designs.
\newblock \emph{European Journal of Operational Research}, 301\penalty0
  (2):\penalty0 576--590, 2022.

\bibitem[Mrkela and Stanimirović(2022)]{Mrkela2022}
L.~Mrkela and Z.~Stanimirović.
\newblock A variable neighborhood search for the budget-constrained maximal
  covering location problem with customer preference ordering.
\newblock \emph{Operational Research}, 22\penalty0 (5):\penalty0 5913 – 5951,
  2022.

\bibitem[Naji-Azimi et~al.(2012)Naji-Azimi, Renaud, Ruiz, and Salari]{Naji2012}
Z.~Naji-Azimi, J.~Renaud, A.~Ruiz, and M.~Salari.
\newblock A covering tour approach to the location of satellite distribution
  centers to supply humanitarian aid.
\newblock \emph{European Journal of Operational Research}, 222\penalty0
  (3):\penalty0 596--605, 2012.

\bibitem[Ortiz-Astorquiza et~al.(2018)Ortiz-Astorquiza, Contreras, and
  Laporte]{ORTIZASTORQUIZA2018791}
C.~Ortiz-Astorquiza, I.~Contreras, and G.~Laporte.
\newblock Multi-level facility location problems.
\newblock \emph{European Journal of Operational Research}, 267\penalty0
  (3):\penalty0 791--805, 2018.

\bibitem[Paul et~al.(2017)Paul, Lunday, and Nurre]{paul2017multiobjective}
N.~R. Paul, B.~J. Lunday, and S.~G. Nurre.
\newblock A multiobjective, maximal conditional covering location problem
  applied to the relocation of hierarchical emergency response facilities.
\newblock \emph{Omega}, 66:\penalty0 147--158, 2017.

\bibitem[Pelegrín and Xu(2023)]{PelXu23}
M.~Pelegrín and L.~Xu.
\newblock Continuous covering on networks: Improved mixed integer programming
  formulations.
\newblock \emph{Omega (United Kingdom)}, 117, 2023.

\bibitem[Toregas et~al.(1971)Toregas, Swain, ReVelle, and
  Bergman]{TorSwaRevBer71}
C.~Toregas, R.~Swain, C.~ReVelle, and L.~Bergman.
\newblock The location of emergency service facilities.
\newblock \emph{Operations Research}, 19\penalty0 (6):\penalty0 1363--1373,
  1971.

\bibitem[Vatsa and Jayaswal(2021)]{VatJay21}
A.~K. Vatsa and S.~Jayaswal.
\newblock Capacitated multi-period maximal covering location problem with
  server uncertainty.
\newblock \emph{European Journal of Operational Research}, 289\penalty0
  (3):\penalty0 1107--1126, 2021.

\bibitem[Yıldırım and Soylu(2023)]{Y2023}
B.~Yıldırım and B.~Soylu.
\newblock Relocating emergency service vehicles with multiple coverage and
  critical levels partition.
\newblock \emph{Computers and Industrial Engineering}, 177, 2023.

\end{thebibliography}

\appendix 
\section{Alternative formulation} \label{App:formulation}

In this appendix, an alternative formulation for the problem is proposed.  For this purpose, the following family of variables  is introduced: 
$$u'_{jt}=\begin{cases}
1,&\mbox{if $j$ offers at least $t$ different products, }\\
0,&\mbox{otherwise,}
\end{cases}j\in J,t\in T(j)\setminus\{0\}.$$

 This family of variables replaces the $u$-variables used in formulation \nameM.  
\begin{align}
\mbox{(\nameMprime)} \quad&\max&&\displaystyle \sum_{i\in I}\sum_{s\in S(i)} w_{is}y_{is}\nonumber\\
&\mbox{s.t.}&&\mbox{\labelcref{con:cubri,con:onesource,con:x,con:y_b}}, \nonumber\\
&&&\sum_{m\in M_J(j)}x_{jm}=\sum_{t\in T(j)\setminus\{0\}}u'_{jt},&&\hspace*{-0.7cm}j\in J,\label{con:x_up}\\
&&&u'_{jt}\leq u'_{j,t-1},&&\hspace*{-0.8cm}j\in J,t\in T(j)\setminus\{0,1\}, \label{con:up}\\
&&&\sum_{j\in J}\sum_{t\in T(j)\setminus\{0\}}(c_{jt}-c_{j,t-1}) u'_{jt}\leq B,\label{con:presu2}\\
&&&u'_{jt}\in \{0,1\},&&\hspace*{-0.8cm}j\in J,t\in T(j)\setminus\{0\}.\label{con:updom}
\end{align}
In this alternative formulation, (\nameMprime), constraints \eqref{con:x_up} and \eqref{con:up} determine the number of different products provided by each $j\in J$. In addition, constraint \eqref{con:presu2} ensures that the operating  cost of second-level facilities do not exceed the budget $B$. As presented in constraints \eqref{con:updom}, the new family of variables $u'$ must be integer.

In preliminary computational experiments, it was found that formulation \nameM performed slightly better than formulation \nameMprime. Thus, the paper has focused on the formulation \nameM.


\section{Best bound of the different approaches} 
 In Table \ref{tb:BS} is reported the best bound found by the exact solution method ($\text{BS}^t$) and the best solution found by the three variants of the heuristic algorithm (C, BP, and BW respectively).  Note that the first column depicts the name of the instance. The names of the instances have the following structure HMCLP\_r\_S\_M\_J\_I\_w, where r represents the number of the instance depending on the random seed generator, S the number of first-level facilities, M the number of different products, J the number of second-level facilities, and I the number of clients. Finally, w represents the approach to generate the preferences,  d indicates that they have been generated using the direct case (dir) and i represents the inverse case (inv). 


\begin{center}
\begin{longtable}{lrrrr}
\caption{Best bounds obtained by the best exact solution approach and the heuristics approaches} \label{tb:BS}\\
\toprule
\textbf{Data}  & \textbf{$\text{BS}^t$ } & \textbf{C} & \textbf{BP} & \textbf{BW}\\ 
\hline
\endfirsthead
\multicolumn{4}{c}%
{\tablename\ \thetable\ -- \textit{Continued from previous page}} \\
\hline
\textbf{Data}  & \textbf{$\text{BS}^t$ } & \textbf{C} & \textbf{BP} & \textbf{BW} \\
\hline
\endhead
\hline \multicolumn{4}{r}{\textit{Continued on next page}} \\ \hline
\endfoot
\hline
\endlastfoot
 
    HMCLP\_1\_30\_8\_100\_100\_d & 45735 & 45652 & 45735 & 45735 \\
    HMCLP\_1\_30\_8\_100\_100\_i & 31665 & 31518 & 31640 & 31640 \\
    HMCLP\_1\_30\_8\_100\_250\_d & 77498 & 76057 & 77498 & 77498 \\
    HMCLP\_1\_30\_8\_100\_250\_i & 48533 & 48396 & 48533 & 48364 \\
    HMCLP\_1\_30\_8\_100\_500\_d & 183559 & 182102 & 183085 & 183085 \\
    HMCLP\_1\_30\_8\_100\_500\_i & 190388 & 188387 & 190013 & 190388 \\
    HMCLP\_1\_30\_8\_250\_1250\_d & 342754 & 330954 & 341607 & 341607 \\
    HMCLP\_1\_30\_8\_250\_1250\_i & 417560 & 412601 & 415847 & 415847 \\
    HMCLP\_1\_30\_8\_250\_250\_d & 148029 & 146864 & 147583 & 147583 \\
    HMCLP\_1\_30\_8\_250\_250\_i & 65086 & 64585 & 65086 & 65033 \\
    HMCLP\_1\_30\_8\_250\_625\_d & 181621 & 179103 & 181399 & 181399 \\
    HMCLP\_1\_30\_8\_250\_625\_i & 148841 & 147286 & 148105 & 148754 \\
    HMCLP\_1\_40\_10\_100\_100\_d & 66686 & 66062 & 66686 & 66601 \\
    HMCLP\_1\_40\_10\_100\_100\_i & 33497 & 33497 & 33497 & 33497 \\
    HMCLP\_1\_40\_10\_100\_250\_d & 131778 & 131547 & 131778 & 131683 \\
    HMCLP\_1\_40\_10\_100\_250\_i & 70219 & 69046 & 70219 & 70219 \\
    HMCLP\_1\_40\_10\_100\_500\_d & 160627 & 157762 & 160627 & 160627 \\
    HMCLP\_1\_40\_10\_100\_500\_i & 136206 & 134897 & 136206 & 136206 \\
    HMCLP\_1\_40\_10\_250\_1250\_d & 732981 & 712820 & 726250 & 726250 \\
    HMCLP\_1\_40\_10\_250\_1250\_i & 342850 & 339299 & 341469 & 341511 \\
    HMCLP\_1\_40\_10\_250\_250\_d & 168827 & 168132 & 168661 & 168672 \\
    HMCLP\_1\_40\_10\_250\_250\_i & 100576 & 99748 & 100357 & 100361 \\
    HMCLP\_1\_40\_10\_250\_625\_d & 267665 & 265617 & 267179 & 267179 \\
    HMCLP\_1\_40\_10\_250\_625\_i & 204226 & 203140 & 204245 & 204245 \\
    HMCLP\_2\_30\_8\_100\_100\_d & 35119 & 34812 & 35119 & 35119 \\
    HMCLP\_2\_30\_8\_100\_100\_i & 34160 & 33943 & 34155 & 34160 \\
    HMCLP\_2\_30\_8\_100\_250\_d & 89595 & 87997 & 89595 & 89595 \\
    HMCLP\_2\_30\_8\_100\_250\_i & 57752 & 57079 & 57752 & 57752 \\
    HMCLP\_2\_30\_8\_100\_500\_d & 160931 & 158314 & 160931 & 159765 \\
    HMCLP\_2\_30\_8\_100\_500\_i & 102458 & 102259 & 101814 & 102458 \\
    HMCLP\_2\_30\_8\_250\_1250\_d & 524139 & 511984 & 523078 & 523078 \\
    HMCLP\_2\_30\_8\_250\_1250\_i & 561359 & 552729 & 558980 & 559289 \\
    HMCLP\_2\_30\_8\_250\_250\_d & 93989 & 93429 & 93830 & 93980 \\
    HMCLP\_2\_30\_8\_250\_250\_i & 74061 & 73784 & 73956 & 73956 \\
    HMCLP\_2\_30\_8\_250\_625\_d & 170850 & 168854 & 169696 & 169938 \\
    HMCLP\_2\_30\_8\_250\_625\_i & 171338 & 169996 & 170601 & 170652 \\
    HMCLP\_2\_40\_10\_100\_100\_d & 51203 & 50597 & 51203 & 51203 \\
    HMCLP\_2\_40\_10\_100\_100\_i & 39675 & 39390 & 39516 & 39675 \\
    HMCLP\_2\_40\_10\_100\_250\_d & 88493 & 87735 & 88493 & 88493 \\
    HMCLP\_2\_40\_10\_100\_250\_i & 56769 & 56546 & 56571 & 56706 \\
    HMCLP\_2\_40\_10\_100\_500\_d & 413896 & 406576 & 412490 & 412490 \\
    HMCLP\_2\_40\_10\_100\_500\_i & 275211 & 274003 & 274561 & 274561 \\
    HMCLP\_2\_40\_10\_250\_1250\_d & 408694 & 390503 & 392090 & 392090 \\
    HMCLP\_2\_40\_10\_250\_1250\_i & 472831 & 471409 & 472704 & 472704 \\
    HMCLP\_2\_40\_10\_250\_250\_d & 119718 & 119097 & 119383 & 119383 \\
    HMCLP\_2\_40\_10\_250\_250\_i & 128496 & 128205 & 128451 & 128451 \\
    HMCLP\_2\_40\_10\_250\_625\_d & 155551 & 154705 & 154763 & 155008 \\
    HMCLP\_2\_40\_10\_250\_625\_i & 204330 & 204057 & 203851 & 203909 \\
    HMCLP\_3\_30\_8\_100\_100\_d & 42821 & 42560 & 42736 & 42736 \\
    HMCLP\_3\_30\_8\_100\_100\_i & 32042 & 31925 & 32042 & 32042 \\
    HMCLP\_3\_30\_8\_100\_250\_d & 46716 & 46584 & 46716 & 46716 \\
    HMCLP\_3\_30\_8\_100\_250\_i & 80076 & 78621 & 79924 & 79951 \\
    HMCLP\_3\_30\_8\_100\_500\_d & 197591 & 194911 & 196319 & 197489 \\
    HMCLP\_3\_30\_8\_100\_500\_i & 143760 & 141964 & 143375 & 143600 \\
    HMCLP\_3\_30\_8\_250\_1250\_d & 461795 & 455939 & 459689 & 459689 \\
    HMCLP\_3\_30\_8\_250\_1250\_i & 382648 & 378286 & 377952 & 379174 \\
    HMCLP\_3\_30\_8\_250\_250\_d & 90643 & 89927 & 90637 & 90637 \\
    HMCLP\_3\_30\_8\_250\_250\_i & 78417 & 78147 & 78348 & 78388 \\
    HMCLP\_3\_30\_8\_250\_625\_d & 197980 & 195675 & 197370 & 197394 \\
    HMCLP\_3\_30\_8\_250\_625\_i & 219897 & 218960 & 219268 & 219724 \\
    HMCLP\_3\_40\_10\_100\_100\_d & 36941 & 36667 & 36941 & 36941 \\
    HMCLP\_3\_40\_10\_100\_100\_i & 35907 & 35522 & 35907 & 35907 \\
    HMCLP\_3\_40\_10\_100\_250\_d & 83286 & 83001 & 83286 & 83087 \\
    HMCLP\_3\_40\_10\_100\_250\_i & 85323 & 83884 & 85323 & 85323 \\
    HMCLP\_3\_40\_10\_100\_500\_d & 221672 & 221303 & 221672 & 221672 \\
    HMCLP\_3\_40\_10\_100\_500\_i & 196447 & 195256 & 196447 & 196447 \\
    HMCLP\_3\_40\_10\_250\_1250\_d & 671084 & 659000 & 670703 & 670526 \\
    HMCLP\_3\_40\_10\_250\_1250\_i & 518092 & 506082 & 517187 & 517187 \\
    HMCLP\_3\_40\_10\_250\_250\_d & 60732 & 60431 & 60718 & 60718 \\
    HMCLP\_3\_40\_10\_250\_250\_i & 76125 & 76085 & 76125 & 76046 \\
    HMCLP\_3\_40\_10\_250\_625\_d & 418415 & 414426 & 417985 & 417985 \\
    HMCLP\_3\_40\_10\_250\_625\_i & 239488 & 235639 & 239133 & 239249 \\
    HMCLP\_4\_30\_8\_100\_100\_d & 28194 & 28033 & 28006 & 28111 \\
    HMCLP\_4\_30\_8\_100\_100\_i & 30742 & 30605 & 30742 & 30742 \\
    HMCLP\_4\_30\_8\_100\_250\_d & 101742 & 100309 & 101742 & 101742 \\
    HMCLP\_4\_30\_8\_100\_250\_i & 78081 & 77865 & 77738 & 78081 \\
    HMCLP\_4\_30\_8\_100\_500\_d & 156226 & 155790 & 156102 & 156226 \\
    HMCLP\_4\_30\_8\_100\_500\_i & 156989 & 156871 & 156982 & 156982 \\
    HMCLP\_4\_30\_8\_250\_1250\_d & 395315 & 387069 & 394923 & 394923 \\
    HMCLP\_4\_30\_8\_250\_1250\_i & 431824 & 426743 & 430421 & 430093 \\
    HMCLP\_4\_30\_8\_250\_250\_d & 80294 & 79767 & 80294 & 80294 \\
    HMCLP\_4\_30\_8\_250\_250\_i & 74495 & 73855 & 74318 & 74318 \\
    HMCLP\_4\_30\_8\_250\_625\_d & 288916 & 283620 & 288140 & 288671 \\
    HMCLP\_4\_30\_8\_250\_625\_i & 144739 & 143576 & 144739 & 144613 \\
    HMCLP\_4\_40\_10\_100\_100\_d & 39382 & 39134 & 39382 & 39382 \\
    HMCLP\_4\_40\_10\_100\_100\_i & 29249 & 29127 & 29249 & 29249 \\
    HMCLP\_4\_40\_10\_100\_250\_d & 53460 & 52987 & 53460 & 53152 \\
    HMCLP\_4\_40\_10\_100\_250\_i & 77273 & 76640 & 77159 & 77159 \\
    HMCLP\_4\_40\_10\_100\_500\_d & 254358 & 251727 & 253823 & 254358 \\
    HMCLP\_4\_40\_10\_100\_500\_i & 199783 & 198209 & 199690 & 199783 \\
    HMCLP\_4\_40\_10\_250\_1250\_d & 407712 & 403415 & 407769 & 407769 \\
    HMCLP\_4\_40\_10\_250\_1250\_i & 473714 & 466112 & 469087 & 469940 \\
    HMCLP\_4\_40\_10\_250\_250\_d & 147578 & 146968 & 147552 & 147604 \\
    HMCLP\_4\_40\_10\_250\_250\_i & 72776 & 72598 & 72776 & 72776 \\
    HMCLP\_4\_40\_10\_250\_625\_d & 169076 & 167395 & 168788 & 168416 \\
    HMCLP\_4\_40\_10\_250\_625\_i & 239942 & 238543 & 239711 & 239711 \\
    HMCLP\_5\_30\_8\_100\_100\_d & 30697 & 30604 & 30615 & 30615 \\
    HMCLP\_5\_30\_8\_100\_100\_i & 22071 & 21790 & 22071 & 22071 \\
    HMCLP\_5\_30\_8\_100\_250\_d & 94955 & 92046 & 94705 & 94955 \\
    HMCLP\_5\_30\_8\_100\_250\_i & 77399 & 76872 & 77278 & 77095 \\
    HMCLP\_5\_30\_8\_100\_500\_d & 232912 & 231269 & 232912 & 232444 \\
    HMCLP\_5\_30\_8\_100\_500\_i & 117214 & 116668 & 117022 & 117045 \\
    HMCLP\_5\_30\_8\_250\_1250\_d & 689224 & 682868 & 687541 & 687541 \\
    HMCLP\_5\_30\_8\_250\_1250\_i & 314315 & 310600 & 313091 & 312804 \\
    HMCLP\_5\_30\_8\_250\_250\_d & 62015 & 61903 & 61835 & 61835 \\
    HMCLP\_5\_30\_8\_250\_250\_i & 79685 & 79346 & 79497 & 79489 \\
    HMCLP\_5\_30\_8\_250\_625\_d & 384674 & 379843 & 382991 & 383223 \\
    HMCLP\_5\_30\_8\_250\_625\_i & 157094 & 156080 & 156373 & 156515 \\
    HMCLP\_5\_40\_10\_100\_100\_d & 52484 & 52484 & 52484 & 52484 \\
    HMCLP\_5\_40\_10\_100\_100\_i & 36898 & 36898 & 36898 & 36898 \\
    HMCLP\_5\_40\_10\_100\_250\_d & 121453 & 120431 & 121338 & 121338 \\
    HMCLP\_5\_40\_10\_100\_250\_i & 75188 & 74527 & 75082 & 75082 \\
    HMCLP\_5\_40\_10\_100\_500\_d & 199054 & 196266 & 199054 & 199007 \\
    HMCLP\_5\_40\_10\_100\_500\_i & 158531 & 157972 & 158531 & 158531 \\
    HMCLP\_5\_40\_10\_250\_1250\_d & 565686 & 550115 & 564685 & 564685 \\
    HMCLP\_5\_40\_10\_250\_1250\_i & 431443 & 426110 & 429715 & 430442 \\
    HMCLP\_5\_40\_10\_250\_250\_d & 123503 & 122750 & 123067 & 123125 \\
    HMCLP\_5\_40\_10\_250\_250\_i & 85829 & 85329 & 85829 & 85791 \\
    HMCLP\_5\_40\_10\_250\_625\_d & 240740 & 238219 & 239871 & 239871 \\
    HMCLP\_5\_40\_10\_250\_625\_i & 218592 & 217085 & 218503 & 218503 \\
\end{longtable}
\end{center}



\end{document}